\newif\ifijnme 
\ijnmefalse 

\ifijnme
\documentclass[times,doublespace]{nmeauth}
\else
\documentclass[preprint,12pt,times]{elsarticle}
\biboptions{sort&compress,comma,square}
\usepackage{amsmath}
\fi

\usepackage{todonotes}
\usepackage{bm}
\usepackage{xcolor}
\usepackage{listings}
\usepackage[linesnumbered,lined,commentsnumbered,ruled]{algorithm2e}

\usepackage[format=plain,indention=.5cm]{caption}
\usepackage[labelformat=simple]{subcaption}

\newcommand*{\gz}[1]{\boldsymbol{#1}}
\newcommand*{\Grad}{\mathrm{Grad}}
\newcommand*{\grad}{\mathrm{grad}}
\renewcommand*{\d}{\mathrm{d}}
\newcommand*{\D}{\mathrm{D}}
\newcommand*{\mcl}[1]{\mathcal{#1}}
\DeclareMathOperator{\trace}{tr}

\begin{document}

\ifijnme
\runningheads{D. Davydov et al.}{A matrix-free approach for finite-strain hyperelastic problems}
\else
\journal{CMAME}
\begin{frontmatter}
\fi

\title{
  A matrix-free approach for finite-strain hyperelastic problems using geometric multigrid
  }

\ifijnme
  \author{
    Denis Davydov\affil{1} \corrauth ,
    Jean-Paul Pelteret\affil{1},
    Daniel Arndt\affil{2},
    Paul Steinmann\affil{1,3}
  }

  \address{
    \affilnum{1}Chair of Applied Mechanics,
    Friedrich-Alexander-Universit\"{a}t Erlangen-N\"{u}rnberg,
    Egerlandstr.\ 5, 91058 Erlangen, Germany \break
    \affilnum{2}Interdisciplinary Center for Scientific Computing (IWR),
    Heidelberg University,
    Im Neuenheimer Feld 205, 69120 Heidelberg, Germany \break
    \affilnum{3}Glasgow Computational Engineering Center (GCEC),
        University of Glasgow, G12 8QQ Glasgow, United Kingdom \break
    }

  \corraddr{E-mail: denis.davydov@fau.de}

\else
  \author[a]{Denis Davydov\corref{cor}}
  \ead{denis.davydov@fau.de}

  \author[a]{Jean-Paul Pelteret}
  \ead{jean-paul.pelteret@fau.de}

  \author[b]{Daniel Arndt}
  \ead{daniel.arndt@iwr.uni-heidelberg.de}

  \author[a,c]{Paul Steinmann}
  \ead{paul.steinmann@fau.de}

  \cortext[cor]{Corresponding author.}

  \address[a]{Chair of Applied Mechanics,
  Friedrich-Alexander-Universit\"{a}t Erlangen-N\"{u}rnberg,
  Egerlandstr.\ 5, 91058 Erlangen, Germany}

  \address[b]{Interdisciplinary Center for Scientific Computing (IWR),
      Heidelberg University,
      Im Neuenheimer Feld 205,
      69120 Heidelberg,
      Germany}

  \address[c]{Glasgow Computational Engineering Center (GCEC),
      University of Glasgow, G12 8QQ Glasgow, United Kingdom
  }
\fi

  \begin{abstract}
    The performance of finite element solvers on modern computer architectures is typically memory bound for sufficiently large problems.
    The main cause for this is that loading matrix elements from RAM into CPU cache is significantly slower than performing the arithmetic operations when solving the problem.
    In order to improve the performance of iterative solvers within the high-performance computing context, so-called matrix-free methods are
    widely adopted in the fluid mechanics community, where matrix-vector products are computed on-the-fly.

    To date, there have been few (if any) assessments into the applicability of the matrix-free approach to problems in solid mechanics.
    In this work, we perform an initial investigation on the application of the matrix-free approach to problems in quasi-static finite-strain hyperelasticity to determine whether it is viable for further extension.
    Specifically, we study different numerical implementations of the finite element tangent operator, and determine whether generalized methods of incorporating complex constitutive behavior might be feasible.
    In order to improve the convergence behavior of iterative solvers, we also propose a method by which to construct level tangent operators
    and employ them to define a geometric multigrid preconditioner.
    The performance of the matrix-free operator and the geometric multigrid preconditioner is compared to the matrix-based implementation with an algebraic multigrid preconditioner on a single node for a representative numerical example of a heterogeneous hyperelastic material in two and three dimensions.
    We conclude that the application of matrix-free methods to finite-strain solid mechanics is promising, and that is it possible to develop numerically efficient implementations that are independent of the hyperelastic constitutive law.
  \end{abstract}

\ifijnme


  \maketitle

\else

\end{frontmatter}
\fi

\section{Introduction}

The performance of finite element solvers on modern computer architectures is typically memory bound for sufficiently large problems.
The main cause for this is that loading pre-computed matrix elements from RAM into CPU cache is significantly slower than performing the arithmetic operations when solving the problem.
In order to improve the performance of iterative solvers, so-called matrix-free methods,
where matrix-vector products are formed on-the-fly \cite{kronbichler12,May2015, Krank2017, Brown2010, Gmeiner2016},
are widely adopted in the fluid mechanics community.
Such methods can also be efficiently implemented on GPUs \cite{Abdelfattah2016, ljungkvist2017multigrid}.

To the best of our knowledge, matrix-free methods are not widely adopted within the solid mechanics community.
We are only aware of a single paper currently in preparation that deals with small strain linear elasticity \cite{Clevenger2018}.
One possible reason is that the tangent operators\footnote{Here and
below 'operator' denotes a finite-dimensional linear operator defined on the vector space $\mathbb R^N$, where $N$ is the number of unknown degrees of freedom in the FE discretization.
Although any linear mapping on a finite-dimensional space is representable by a matrix, we adopt the `operator' term to avoid confusion between its matrix-free and matrix-based numerical implementation.
} (stemming from linearization of the nonlinear balance equations) are
more elaborate as compared to the operators used often in both linear and non-linear fluid mechanics, and it is not obvious whether matrix-free methods can be advantageous in this case.
Another reason could be the frequent use of lower-order elements in solid mechanics due to lacking smoothness of the underlying solution.
The finite element method (FEM) with linear or quadratic elements (that reduce the potential for locking), or mixed or enhanced formulations (that avoid locking) are often adopted in the solid mechanics community,
whereas higher order elements are rarely employed.
This can potentially reduce the advantage of matrix-free methods using sum-factorization techniques, that are generally more competitive for higher order elements \cite{kronbichler12,kronbichler2017fast,muthing2017high}.
However, we note that, even with linear FEM, large scale computations with $10^{12}$ unknowns inevitably require a matrix-free approach as there is simply not enough memory to store the sparse tangent matrix \cite{Gmeiner2016}. Therefore, we believe that for large scale computations in solid mechanics it is crucial to consider matrix-free approaches and develop efficient multigrid solvers.

In this work, we perform a preliminary investigation of different implementations of the finite-strain hyperelastic tangent operator with respect to their computational cost
and performance on a single node of a high performance computing cluster.
The \texttt{deal.II} \cite{dealII90} finite element library, which provides MPI parallelization together with a generalized matrix-free framework that incorporates SIMD vectorization and an interface to high-performance, MPI distributed matrix-based libraries, is used for this study.
In order to improve the convergence of iterative solvers, we also propose a method by which to construct level matrix-free tangent operators
and employ them to define a \mbox{geometric} multigrid preconditioner.

The paper is organized as follows: In Section \ref{sec:theory} we briefly introduce the partial differential equations used in  finite-strain solid mechanics.
Sections \ref{sec:fe} and \ref{sec: Description of numerical framework} cover the finite element discretization and provide details of the numerical framework with which this study is implemented.
In \ref{sec:mf} we describe the different matrix-free operator implementations that are evaluated in this study, and in section \ref{sec:gmg} we propose a geometric multigrid preconditioner that is suitable for the matrix-free implementation of finite-strain hyperelasticity.
The performance of the matrix-free operator and the geometric multigrid preconditioner is compared to the matrix-based implementation with an algebraic multigrid preconditioner for a representative numerical example of a heterogeneous hyperelastic material simulated in two and three dimensions in Section \ref{sec:example}.
Lastly, the results are summarized in Section \ref{sec:summary}.

\section{Theoretical background}
\label{sec:theory}

\subsection{Weak form}
The deformation of a body $\mcl B$ from the referential configuration $\mcl B_0$ to the spatial configuration $\mcl B_t$ at time $t$
is defined via the deformation map $\gz \varphi_t: \mcl B_0 \rightarrow \mcl B_t$, which places material points of $\mcl B$ into the Euclidean space $\mathbb E^3$.
The spatial location of a material particle $\gz X$ is given by $\gz x = \gz \varphi_t (\gz X)$.
The displacement field reads $\gz u = \gz x - \gz X$.

The tangent map is linear such that
$\d \gz x = \gz F \cdot \d \gz X$,
where $\gz F := \Grad \, \gz \varphi \equiv \gz I + \Grad \, \gz u$ is called the deformation gradient and $\gz I$ is the second-order unit tensor.
The mapping has to be one-to-one and must exclude self-penetration. Consequently, the Jacobian $J = \det \gz F > 0$ has to be positive.
We shall denote the push-forward transformation of rank-2 and rank-4 tensors $\mathbf{A}$ and $\boldsymbol{\mathcal{A}}$ by a linear transformation $\chi$ defined as
\begin{align}
  \chi\left( \mathbf{A} \right)_{ij}
  &= F_{iA} A_{AB} F_{jB} \\
  \chi\left( \boldsymbol{\mathcal{A}} \right)_{ijkl}
  &= F_{iA} F_{jB} \mathcal{A}_{ABCD} F_{kC} F_{lD}.
\end{align}
Note that $\d \gz x = \chi(\d \gz X)$.

In what follows, we parameterize the material behavior using the right Cauchy-Green tensor $\gz C := \gz F^T \cdot \gz F$.
We will also use the Green-Lagrange strain tensor $\gz E:=\frac{1}{2}\left[\gz C - \gz I\right]$ and the left Cauchy-Green tensor $\gz b := \gz F \cdot \gz F^T$.
Clearly $J = \sqrt{\det \gz C}$ and $\partial (\bullet)/\partial \gz E = 2\, \partial (\bullet)/\partial \gz C$.
For conservative systems without body forces, the total potential energy functional $\mcl E$ is introduced as
\begin{align}
\mcl E =
\int_{\mcl B_0} \mcl \psi(\gz F) \, \d V \,
- \int_{\partial \mcl B_0^N} \overline{\gz T} \cdot \gz u \, \d S \, ,
\end{align}
where $\overline{\gz T}$ is the prescribed loading at the Neumann part of the boundary $\partial \mcl B_0^N$ in the referential configuration, $\mcl \psi$ denotes the strain-energy per unit reference volume and $\gz u$ should satisfy the prescribed Dirichlet boundary conditions $\gz u = \overline{\gz u}$ on $\partial \mcl B_0^D := \partial \mcl B_0 \setminus \partial \mcl B_0^N \neq \emptyset $.

The principle of stationary potential energy at equilibrium requires that the directional derivative with respect to the displacement
\begin{align}
\D_{\displaystyle \delta \gz u} \mcl E :=
\frac{\d}{\d \epsilon} \mcl E (\gz u + \epsilon \delta \gz u) \Bigr\rvert_{\epsilon=0}  = 0 \, \qquad \forall \delta \gz u \, .
\label{eq:stationary}
\end{align}
vanishes for all directions $\delta \gz u$ which satisfy homogeneous Dirichlet boundary conditions.
This leads to the following scalar-valued non-linear equation
\begin{align}
F(\gz u, \delta \gz u) =
\int_{\mcl B_0} \gz P : \Grad \, \delta \gz u \, \d V \,
-
\int_{\partial \mcl B_0^N} \overline{\gz T} \cdot \delta \gz u \, \d S
= 0 \, ,
\end{align}
where $\gz P := \partial \mcl \psi / \partial \gz F$ is the Piola stress tensor.
The double contraction in the first term can be re-written in terms of the symmetric Kirchhoff stress tensor $\gz \tau := \gz P \cdot \gz F^T$ as
\begin{align}
\gz P : \Grad \, \delta \gz u =
\left[\gz \tau \cdot \gz F^{-T}\right] : \Grad \, \delta \gz u =
\gz \tau : \left[\Grad \, \delta \gz u \cdot \gz F^{-1}\right] =
\gz \tau : \grad \, \delta \gz u \, ,
\end{align}
and therefore
\begin{align}
  F(\gz u, \delta \gz u) =
  \int_{\mcl B_0} \gz \tau : \grad^{s} \, \delta \gz u \, \d V \,
  -
  \int_{\partial \mcl B_0^N} \overline{\gz T} \cdot \delta \gz u \, \d S
  = 0 \, ,
\label{eq:weak_form}
\end{align}
where $\grad^{s} \, \delta \gz u := \frac{1}{2}\left[ \grad \, \delta \gz u + ( \grad \, \delta \gz u)^T\right]$
is involved due to the symmetry of $\gz \tau$.
Note that $\gz \tau$ is the push-forward transformation of the Piola-Kirchhoff stress $\gz S := \partial \mcl \psi / \partial \gz E \equiv 2 \partial \mcl \psi / \partial \gz C$,
that is $\gz \tau = \chi(\gz S) = \gz F \cdot \gz S \cdot \gz F^T$.

\subsection{Linearization}

In order to solve \eqref{eq:weak_form} using Newton's method, a first order approximation around a given solution field $\overline{\gz u}$ is required such that
\begin{align}
F(\overline{\gz u} + \Delta \gz u, \delta \gz u) \approx
F(\overline{\gz u}, \delta \gz u) + \D_{\displaystyle \Delta \gz u} F(\overline{\gz u}, \delta \gz u) ,
\end{align}
where $\D_{\displaystyle \Delta \gz u}(\bullet)$ denotes the directional derivative in the direction $\Delta \gz u$.
For conservative traction boundary conditions, the directional derivative is given by
\begin{equation}
\begin{split}
\D_{\displaystyle \Delta \gz u} F(\overline{\gz u}, \delta \gz u)
&=
\int_{\mcl B_0}
\D_{\displaystyle \Delta \gz u} \left(\gz F \cdot \gz S \cdot \gz F^T\right)  :
\overline{\grad^s} \, \delta \gz u
\, \d V
\\
& +
\int_{\mcl B_0}
\overline{\gz \tau} :
\left[
  \Grad \, \delta \gz u \cdot
  \D_{\displaystyle \Delta \gz u} \gz F^{-1}
\right] \d V.
\end{split}
\label{eq:tangent_pre}
\end{equation}
Following \cite{Wriggers2008}, this can be simplified to\footnote{To that end,
$\!\D_{\displaystyle \Delta \gz u} \gz E = \frac{1}{2}\left[\D_{\displaystyle \Delta \gz u} \gz F^T \cdot \gz F + \gz F^T \cdot \D_{\displaystyle \Delta \gz u} \gz F\right]$, $\D_{\displaystyle \Delta \gz u} \gz F = \Grad \Delta \gz u\;$  and
$\!\D_{\displaystyle \Delta \gz u} \gz F^{-1} = - \gz F^{-1} \cdot \D_{\displaystyle \Delta \gz u} \gz F \cdot \gz F^{-1}$
are employed together with
$\D_{\displaystyle \Delta \gz u} \gz S = 2 \partial \gz S / \partial \gz C : \D_{\displaystyle \Delta \gz u} \gz E$.
}
\begin{equation}
  \begin{split}
\D_{\displaystyle \Delta \gz u} F(\overline{\gz u}, \delta \gz u)
  &=
  \int_{\mcl B_0} \overline{\grad^s} \Delta \gz u : J \boldsymbol{\mathcal{C}} : \overline{\grad^s} \, \delta \gz u \, \d V \\
  &+
  \int_{\mcl B_0}
  \overline{\grad}\delta \gz u :
  \left[
  \overline{\grad} \Delta \gz u \cdot
  \overline{\gz \tau}
  \right]
  \d V.
\end{split}
\label{eq:tangent}
\end{equation}
Here $\overline{(\bullet)}$ is used to denote quantities evaluated using the displacement field $\overline {\gz u}$, and the fourth-order material part of the spatial tangent stiffness tensor is the push forward of the material part of the referential tangent stiffness tensor $J \boldsymbol{\mathcal{C}} = \chi\left( 4 \frac{d^{2} \psi \left( \mathbf{C} \right)}{d \mathbf{C} \otimes d \mathbf{C}} \right)$.
The first term in \eqref{eq:tangent} is related to the linearization of the Piola-Kirchhoff stress $\gz S$ and is therefore called material part of the directional derivative.
The second term in \eqref{eq:tangent} is called the geometric part of the directional derivative since it originates from the linearization of $\gz F$, $\gz F^T$ and $\gz F^{-1}$.

\subsection{Constitutive modelling
\label{sec: Constitutive modelling}
}

For the current study, we use the compressible Neo-Hookean model
\begin{gather}
\psi \left( \mathbf{C} \right)
  = \frac{\mu}{2} \left[ \trace{\mathbf{C}} - \trace{\mathbf{I}} - 2 \ln\left( J \right) \right]
  + \lambda \ln^{2}\left( J \right)
\end{gather}
where $\mu$ and $\lambda$ denote the shear modulus and Lam\'{e} parameter respectively.
Derived using statistical thermodynamics applied to cross-linked polymers, the Neo-Hookean model \cite{Treloar1975a,Treloar1976a} is commonly used to model rubber-like materials in the finite-strain regime.
It can be shown that the first and second derivatives of the strain energy function are given by
\begin{align}
\frac{d \psi \left( \mathbf{C} \right)}{d \mathbf{C}}
  &= \frac{\mu}{2} \mathbf{I} - \frac{1}{2} \left[ \mu - 2\lambda\ln\left( J \right) \right] \mathbf{C}^{-1} \\
\frac{d^{2} \psi \left( \mathbf{C} \right)}{d \mathbf{C} \otimes d \mathbf{C}}
  &= \frac{1}{2}\left[ \mu - 2\lambda\ln\left( J \right) \right] \left[ - \frac{d \mathbf{C}^{-1}}{d \mathbf{C}} \right]
  + \frac{\lambda}{2} \mathbf{C}^{-1} \otimes \mathbf{C}^{-1}
\end{align}
The Kirchhoff stress and its associated fourth-order material part of the spatial tangent tensor are
\begin{gather}
\boldsymbol{\tau}
  \equiv J \boldsymbol{\sigma}
  = \chi\left( 2 \frac{d \psi \left( \mathbf{C} \right)}{d \mathbf{C}} \right)
  = \mu \mathbf{b} - \left[ \mu - 2\lambda\ln\left( J \right) \right] \mathbf{I}
\end{gather}
and
\begin{gather}
J \boldsymbol{\mathcal{C}}
  = \chi\left( 4 \frac{d^{2} \psi \left( \mathbf{C} \right)}{d \mathbf{C} \otimes d \mathbf{C}} \right)
  = 2 \left[ \mu - 2\lambda\ln\left( J \right) \right] \boldsymbol{\mathcal{S}}
  + 2 \lambda \mathbf{I} \otimes \mathbf{I}
\end{gather}
where $\boldsymbol{\mathcal{S}}$ is the fourth-order symmetric identity tensor
with $S_{ijkl}=\frac{1}{2}\left[\delta_{ik}\delta_{jl}+\delta_{il}\delta_{jk}\right]$.
The action that $J \boldsymbol{\mathcal{C}}$ performs when contracted with an arbitrary rank-2 symmetric tensor is therefore
\begin{gather}
J \boldsymbol{\mathcal{C}} : \left( \bullet \right)
  = 2 \left[ \mu - 2\lambda\ln\left( J \right) \right] \left( \bullet \right)
  + 2 \lambda \trace\left( \bullet \right) \mathbf{I}.
\label{eq:simplified_action}
\end{gather}

\section{Finite Element discretization}
\label{sec:fe}

We now introduce a FE triangulation $\mathcal{B}^h_0$ of $\mcl B_0$ and
the associated FE space of continuous piecewise elements with polynomial space of fixed degree. 
The displacement fields are given in
a vector space spanned by standard vector--valued FE basis functions $\gz N_i(\gz x)$ (e.g. polynomials with local support on a patch of elements):
\begin{alignat}{2}
       \gz u^h &=:  \sum_{i \in \mcl I} u_i \gz N_i (\gz X) \quad \quad \quad
\delta \gz u^h &&=: \sum_{i \in \mcl I} \delta u_i \gz N_i (\gz X) \,,
\end{alignat}
where the superscript $h$ denotes that this representation is related to the FE mesh with size function $h(\gz X)$ and $\mcl I$ is the set of unknown degrees of freedom (DoF).

The Newton-Raphson solution approach is adopted for the nonlinear problem considered here.
Given the current trial solution field $\overline{\gz u^h}$, the correction $\Delta \gz u^h$ field is obtained as the solution to the following system of equations (see \eqref{eq:weak_form} and \eqref{eq:tangent}):
\begin{align}
  \sum_{j \in \mcl I} A_{ij} \Delta u_j &= - F_i  \label{eq:linear_system} \\
  A_{ij} \equiv a(\gz N_i, \gz N_j) &=
  \int_{\mcl B_0}
  \underbrace{
  \overline{\grad^s} \gz N_i : J \boldsymbol{\mathcal{C}} : \overline{\grad^s} \, \gz N_j
  +
  \overline{\grad}\gz N_i :
  \left[
  \overline{\grad} \gz N_j \cdot
  \overline{\gz \tau}
  \right]
  }_{\displaystyle \tilde{a}(\gz N_i, \gz N_j)}
  \d V
  \label{eq:algebraic_tangent}
  \\
  F_i &=
  \int_{\mcl B_0} \overline{\gz \tau} : \grad^{s} \, \gz N_i \, \d V \,
  -
  \int_{\partial \mcl B_0^N} \overline{\gz T} \cdot \gz N_i \, \d S.
\end{align}
Here, $\gz A$ is the discrete tangent operator, $\gz F$ is the discrete gradient of the potential energy
and $\tilde{a}(\gz N_i, \gz N_j)$ is a representation for the integrand in the bilinear form $a(\gz N_i, \gz N_j)$ for the trial solution field $\overline{\gz u^h}$.
Note that $\overline{\grad^s} \gz N_i$ and $\overline{\grad} \gz N_i$ are gradients of shape functions in the spatial configuration, whereas the integration is done in the referential configuration.

\section{Description of the utilized numerical framework
\label{sec: Description of numerical framework}
}

The discretized system of linear equations formulated in Section \ref{sec:fe}, in conjunction with the constitutive model described in Section \ref{sec: Constitutive modelling} were implemented using \texttt{deal.II} version 9.0 \cite{dealII90}, an open-source finite element library.
This library offers a number of paradigms by which to parallelize a finite-element code, and as a member of the xSDK (Extreme-scale Scientific Software Development Kit) ecosystem \citep{Bartlett2017b} aims to support exascale computing.
In this work we employ \texttt{deal.II}'s interface to \texttt{p4est} \cite{p4est} library to perform a standard domain decomposition, whereafter each MPI process ``owns'' only a subset of elements and the mesh is distributed across all MPI processes.
\begin{figure}[!ht]
  \centering
  \begin{subfigure}[b]{0.8\textwidth}
    \centering
    \includegraphics[width=\textwidth]{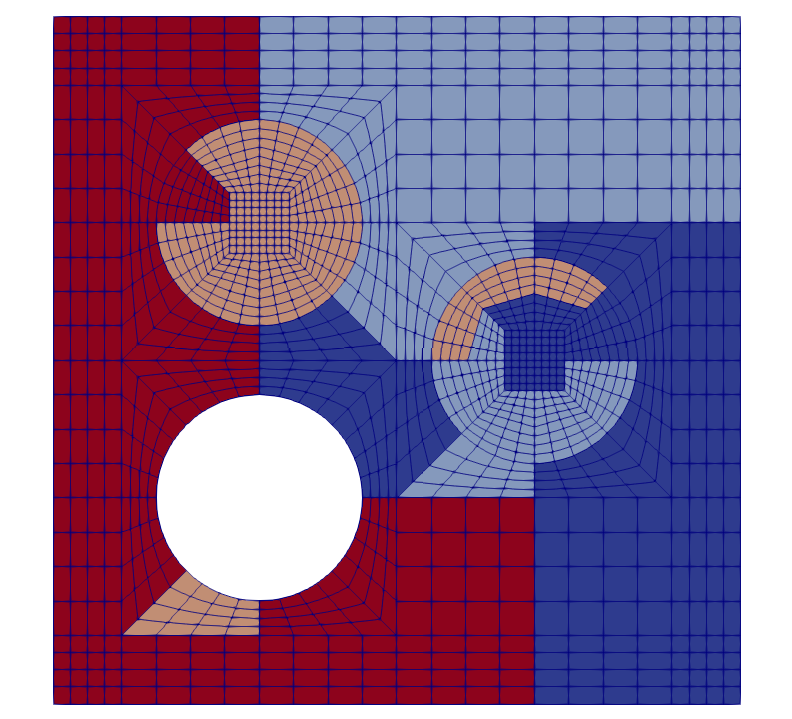}
  \end{subfigure}
  \caption{2D mesh after two global refinements distributed into three MPI processes. The color indicates the process owning the element.}%
  \label{fig:miehe_fine_level}
\end{figure}
Figure \ref{fig:miehe_fine_level} illustrates a distributed MPI decomposition on the fine mesh level for the 2D problem that will be introduced later in Section \ref{sec:example}.

In practice, this implies that parallelization of the assembly operation (for the matrix-based approach) and the linear solver (as well as various pre- and post-processing steps), both of which are the focal points of investigation in this manuscript, has been achieved using MPI.
For the linear solver, \texttt{deal.II}'s implementation of the preconditioned conjugate gradient (CG) solver for symmetric positive-definite systems is consistently used, leaving the choice for method of assembly and linear system preconditioning remaining.

The discrete tangent operator described by \eqref{eq:algebraic_tangent} can be implemented in one of two ways, namely through a classical matrix-based approach or, alternatively, a matrix-free approach.
For the matrix-based approach adopted in this study, the Trilinos \cite{Heroux2005} library was leveraged;
the system tangent matrix is stored in a distributed \textit{Epetra\_FECrsMatrix} from the Epetra package \cite{Heroux2005b} sparse data structure and the linear solver was preconditioned with an algebraic multigrid (AMG) preconditioner from the ML package \cite{Gee2006a}.
This is the most highly performant matrix-based approach using the Trilinos framework that, to date, is possible in  \texttt{deal.II}.

The second implementation, which is discussed in detail in the following Sections, utilized a matrix-free approach in conjunction with a geometric-multigrid preconditioner.
As will be detailed later, the pre-existing matrix-free framework that was used offers further opportunity for low-level parallelism that we exploit.

\section{Matrix-free operator evaluation}
\label{sec:mf}
Classically, in order to solve \eqref{eq:linear_system} the matrix $\gz A$ and the residual force vector $\gz F$ corresponding to the discretization are assembled (that is, the non-zero matrix elements $A_{ij}$ are individually computed and stored) and a direct or iterative solver is used to solve the linear system.
In this case, the matrix-vector product $\gz A \gz x$ results from the product of the individual matrix elements with the elements in the vector, and is usually implemented with the aid of specialized linear algebra packages.
On modern computer architectures however, loading sparse matrix data into the CPU registers is significantly slower than performing the arithmetic operations when solving.
For this reason, recent implementations often focus on so-called matrix-free approaches where the matrix is not assembled but rather the result of its operation on a vector is directly calculated.
The idea is to perform the operations within the solver on-the-fly rather than loading matrix elements from memory.
For iterative solvers, this is possible since it is sufficient to compute matrix-vector products alone.
The matrix-vector product $\gz A \gz x$ (i.e. the action of operator $\gz A$ on a vector $\gz x$) can be expressed by
\begin{align}
  \begin{split}
 (\gz A \gz x)_i &= \sum_j a(\gz N_i,\gz N_j) x_j \\
        &\approx \sum_K \sum_q \sum_j \tilde{a}(\gz N_i,\gz N_j)(\gz \xi_q) x_j w_q J^K_q
  \end{split}
  \label{eq:mf_vmult}
\end{align}
where $\tilde{a}(\gz N_i,\gz N_j)(\gz \xi_q)$ is a representation for the evaluation of the bilinear form in one quadrature point $\gz \xi_q$, $J^K_q$ is the Jacobian of the mapping from the isoparametric element to the element $K$ and
$w_q$ is the quadrature weight.
For the sake of demonstration, let us assume that $\gz A$ represents a mass operator with scalar valued shape functions $\{ N_i(\gz x) \}$. In this case, it holds that
\begin{align*}
 \tilde{a}(N_i,N_j)(\gz \xi_q) = N_i(\gz \xi_q)N_j(\gz \xi_q).
\end{align*}
Further assuming that the number of quadrature points $n$ in one-dimension matches the degree of the ansatz space plus one, the evaluation of all of the shape functions ($n^d$)
in all quadrature points ($n^d$) has complexity $\mathcal{O}(n^{2d})$ in $d$ space dimensions.
For the total matrix-vector product we get $\mathcal{O}(n^{2d})$,
and this is the same as for a regular matrix-vector product where the entries are already precomputed. Assuming tensor product ansatz spaces and a tensor product quadrature rule, evaluation of this operations can be performed more efficiently using a technique called ``sum factorization''.

\begin{figure}[!ht]
  \centering
  \begin{subfigure}[b]{0.45\textwidth}
      \centering
      \includegraphics[width=\textwidth]{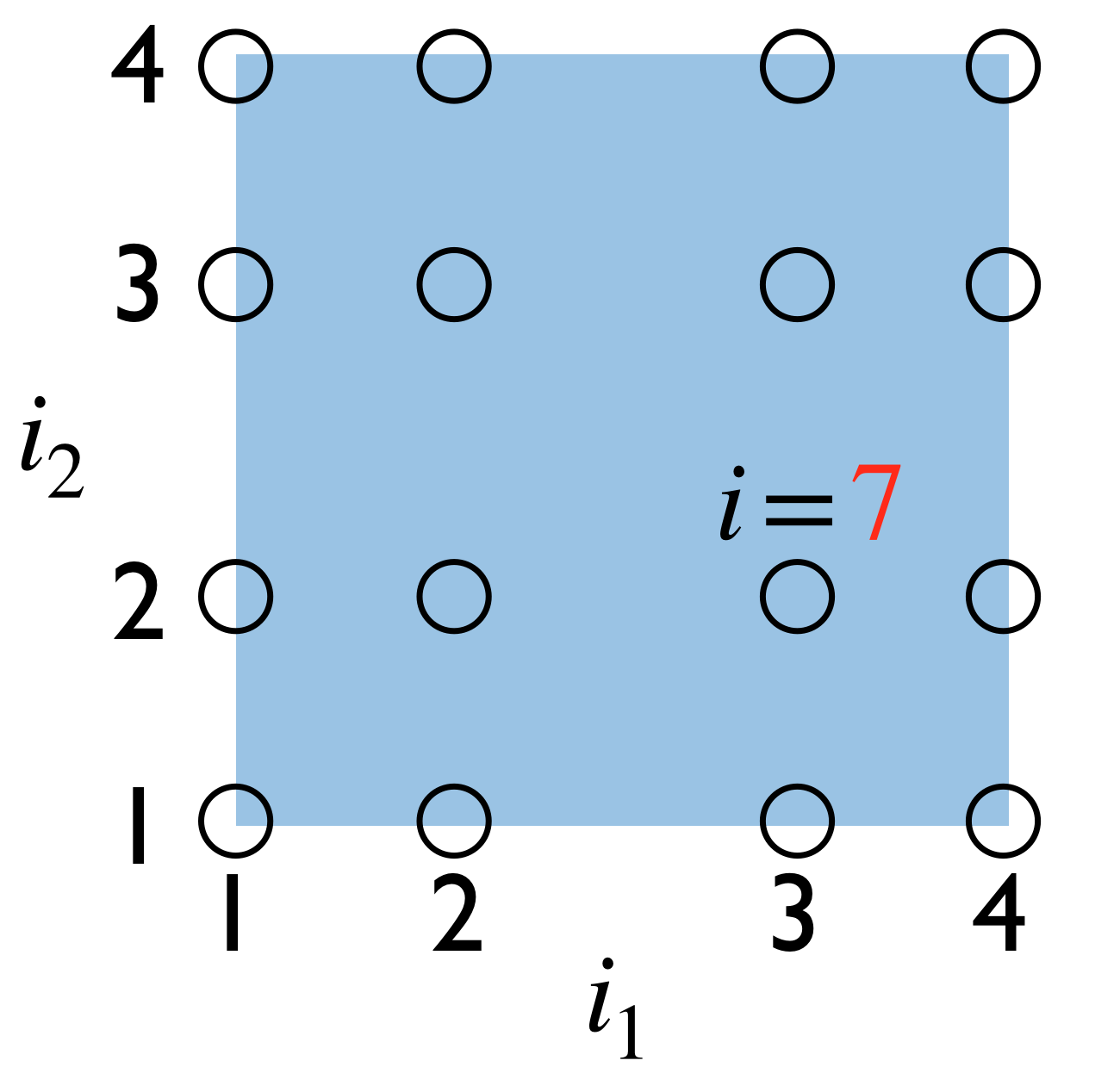}
      \caption{Nodes of cubic Lagrange element}
  \end{subfigure}
  \begin{subfigure}[b]{0.45\textwidth}
    \centering
    \includegraphics[width=\textwidth]{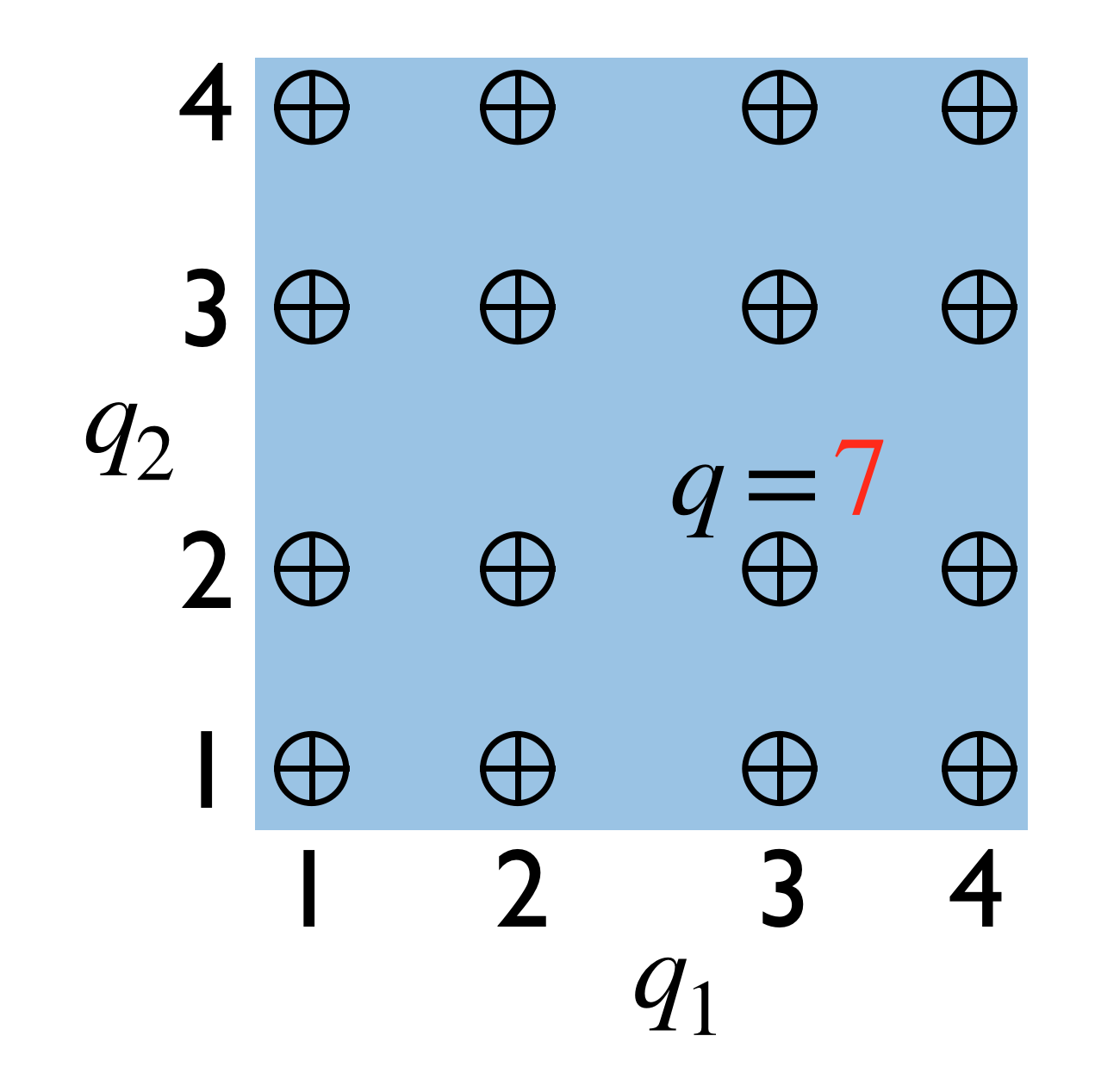}
    \caption{Gaussian quadrature points}
  \end{subfigure}
  \caption{Illustration of tensorial indexing for nodes and quadrature points. In this sketch a node/quadrature point with the number $7$ can be indexed with a multi-index $\{3,2\}$.}%
  \label{fig:multi_index}
\end{figure}

Let us illustrate this in 2D. Tensor-product quadrature points can be expressed as a combination of one dimensional quadrature points
\begin{align}
  \gz \xi_q = (\widetilde{\xi}_{q_1}, \widetilde{\xi}_{q_2}),
\end{align}
where for each quadrature points $q$ we can associate a multi-index $(q_1,q_2)$.
Weights associated with the quadrature formula $w_q$ can be expressed as a product of one-dimensional weights
\begin{align}
  w_q = \widetilde{w}_{q_1} \widetilde{w}_{q_2} \, .
\end{align}
Shape functions can also be expressed as product of one dimensional basis functions
\begin{align}
  N_i(\gz \xi_q) =
  \widetilde{N}_{i_1}(\widetilde{\xi}_{q_1})
  \widetilde{N}_{i_2}(\widetilde{\xi}_{q_2})
\end{align}
where for each DoF index $i$ we can associate a multi-index $(i_1, i_2)$.
See Figure \ref{fig:multi_index} for an illustration.
Therefore, the result of the application of the scalar-valued mass operator, represented as a vector $y_i$, on a single element $K$ can be written as
\begin{align*}
  y_{i} \equiv Y_{i_1\,i_2} & = \sum_{q_1} \sum_{q_2} \sum_{j_1} \sum_{j_2}
  \widetilde{N}_{i_1}(\widetilde{\xi}_{q_1})
  \widetilde{N}_{i_2}(\widetilde{\xi}_{q_2})
  \widetilde{N}_{j_1}(\widetilde{\xi}_{q_1})
  \widetilde{N}_{j_2}(\widetilde{\xi}_{q_2})
  X_{j_1\,j_2}
  \widetilde{w}_{q_1} \widetilde{w}_{q_2} J^K_{q_1\,q_2}
  \\
  &=
  \sum_{q_1} \widetilde{N}_{i_1}(\widetilde{\xi}_{q_1}) \widetilde{w}_{q_1}
  \sum_{q_2} \widetilde{N}_{i_2}(\widetilde{\xi}_{q_2}) \widetilde{w}_{q_2}
  J^K_{q_1\,q_2}
  \left[
    \sum_{j_1}
    \widetilde{N}_{j_1}(\widetilde{\xi}_{q_1})
    \sum_{j_2}
    \widetilde{N}_{j_2}(\widetilde{\xi}_{q_2})
    X_{j_1\,j_2}
  \right]\, \forall i_1 \, i_2
\end{align*}
Note that $x_j \equiv X_{j_1\,j_2}$; the former accesses the element source vector using a linear index $j$, whereas the latter uses the associated multi-index accesses $\{j_1\,j_2\}$.
These four loops all have the same structure.
We either keep the shape function fixed and iterate over the quadrature points in one spatial direction (in the first two sums) or keep the quadrature point fixed and iterate over all the shape functions in one spatial direction (in the last two sums).
Since each of these $2\times d$ loops has a complexity of $n$ and we perform them for $n^d$ values in quadrature points or values in degrees of freedoms simultaneously, this results in an algorithm with an arithmetic complexity of $\mathcal{O}(n^{d+1})$.
Hence, we can expect that such a matrix-free approach is faster than a regular matrix-based approach, especially for the 3D case, and where the discretization and evaluation employs high polynomial degree basis functions with the corresponding quadrature rule.
More information on matrix-free techniques can be found in \cite{kronbichler12,vos10}.

For the performance of the matrix-free operator, it is crucial to find a good intermediate point between precomputing everything (caching) and evaluating all the quantities on-the-fly.
In order to obtain gradients with respect to the spatial configuration,
that appear in the tangent operator \eqref{eq:algebraic_tangent} of the finite-strain elasticity problem derived with a Lagrangian formulation,
we used the \textit{MappingQEulerian} class from the \texttt{deal.II} \cite{dealII90} library.
In addition to the standard mapping from the reference (isoparametric) element to the element in real space, this class computes the mapping to the spatial configuration given a displacement field.
Due to the specifics of matrix-free operator evaluation implemented in \texttt{deal.II}, the integration is eventually performed over the spatial configuration.
Therefore, we have to additionally divide by the Jacobian $J$ of the deformation map at a given quadrature point to express the integral in the referential configuration.

Below, we propose three algorithms to implement the tangent operator of finite-strain elasticity (see \eqref{eq:algebraic_tangent}) using the matrix-free approach, each of which have a different level of abstraction, algorithmic complexity, and memory requirements.
They are introduced in order of lowest memory footprint to largest:

\begin{algorithm}[!ht]
  \SetKwInOut{Input}{Given}
  \SetKwInOut{Output}{Return}
  \Input{Source FE vector $\gz x$, current FE solution $\overline{\gz u^h}$,
  cached $c_1 := \mu - 2 \lambda \log(J)$ for each cell and quadrature point}
  \Output{action of the FE tangent operator \eqref{eq:algebraic_tangent} on $\gz x$}
  \ForEach{ element $K \in \Omega^h$ }{
        extract local vector values on this element $\overline{\gz u^h}_K$, $\gz x_K$ \;
        evaluate the following tensors at each quadrature point using sum factorization :\\
        \enskip $\Grad \, \overline{\gz u^h}_K$ \tcp*{2nd order}
        \enskip $\gz g_s := \overline{\grad^s} \gz x_K$ \tcp*{2nd order symmetric}
        \enskip $\gz g := \overline{\grad} \gz x_K$ \tcp*{2nd order}
        \ForEach{quadrature point $q$ on $K$}{
          evaluate $\gz F = \gz I + \Grad \, \overline{\gz u^h}$ \tcp*{2nd order}
          evaluate $J = \rm{det}(\gz F)$ \tcp*{scalar}
          evaluate $\gz b = \gz F \cdot \gz F^T$ \tcp*{2nd order symmetric}
          evaluate $\gz \tau = \mu \gz b - c_1 \gz I$ \tcp*{2nd order symmetric}
          evaluate $\boldsymbol{\mathcal{G}} \gz g := \gz g \cdot \gz \tau/J$ \tcp*{2nd order}
          queue $\boldsymbol{\mathcal{G}} \gz g$ for contraction $\overline{\grad} \gz N_i : \boldsymbol{\mathcal{G}} \gz g$\;
          evaluate $\boldsymbol{\mathcal{C}}\gz g_s := \left[2 c_1 \gz g_s + 2 \lambda\, \rm{tr}( \gz g_s) \gz I \right]/J$ \tcp*{2nd order symmetric}
          queue $\boldsymbol{\mathcal{C}} \gz g_s$ for contraction $\overline{\grad^s} \gz N_i : \boldsymbol{\mathcal{C}}\gz g_s$ \;
        }
      evaluate queued contractions using sum factorization \;
      distribute results to the destination vector
  }
  \caption{Matrix-free tangent operator: cache scalar quantities only}
  \label{alg:mf_scalar}
\end{algorithm}

Algorithm \ref{alg:mf_scalar} caches a single scalar which involves the logarithm of the Jacobian -- a computationally demanding operation compared to additions and multiplications. As evident from Section \ref{sec:theory}, we need to re-evaluate the Kirchhoff stress at each quadrature point in order to apply the tangent operator. This in turn requires evaluating the deformation gradient $\gz F$ and therefore gradients with respect to the referential configuration.

\begin{algorithm}[!ht]
  \SetKwInOut{Input}{Given}
  \SetKwInOut{Output}{Return}
  \Input{Source FE vector $\gz x$,
  cached
  $c_1 := 2\left[\mu - 2 \lambda \log(J)\right]/J$,
  $c_2 := 2\lambda/J$ and $\gz \tau/J$ for each cell and quadrature point,
  evaluated based on $\overline{\gz u^h}$}
  \Output{action of the FE tangent operator \eqref{eq:algebraic_tangent} on $\gz x$}
  \ForEach{ element $K \in \Omega^h$ }{
        extract local vector values on this element $\gz x_K$ \;
        evaluate the following tensors at each quadrature point using sum factorization :\\
        \enskip$\gz g_s := \overline{\grad^s} \gz x_K$ \tcp*{2nd order symmetric}
        \enskip$\gz g := \overline{\grad} \gz x_K$ \tcp*{2nd order}
        \ForEach{quadrature point $q$ on $K$}{
          evaluate $\boldsymbol{\mathcal{G}}\gz g := \gz g \cdot \left[\gz \tau/J \right]$ \tcp*{2nd order}
          queue $\boldsymbol{\mathcal{G}} \gz g$ for contraction $\overline{\grad} \gz N_i : \boldsymbol{\mathcal{G}} \gz g$ \;
          evaluate $\boldsymbol{\mathcal{C}} \gz g_s := c_1 \gz g_s + c_2 \, \rm{tr}(\gz g_s) \gz I$ \tcp*{2nd order symmetric}
          queue $\boldsymbol{\mathcal{C}} \gz g_s$ contraction $\overline{\grad^s} \gz N_i : \boldsymbol{\mathcal{C}} \gz g_s $ \;
        }
        evaluate queued contractions using sum factorization \;
        distribute results to the destination vector
  }
  \caption{Matrix-free tangent operator: cache second-order Kirchhoff stress $\gz \tau$ and thereby avoid the need to evaluate referential quantities like $\gz F$ at when executed.}
  \label{alg:mf_tensor2}
\end{algorithm}

In order to avoid recalculation of the Kirchhoff stress, Algorithm \ref{alg:mf_tensor2} caches its value for each element and quadrature point. In this case, we also avoid the need to evaluate gradients of the displacement field with respect to the referential configuration.
This algorithm also utilizes the chosen constitutive relationship in that the operation of $J \boldsymbol{\mathcal{C}}$ remains directly expressed using \eqref{eq:simplified_action}.
Therefore, the action of the material part of the fourth-order spatial tangent stiffness tensor $J \boldsymbol{\mathcal{C}}$ on the second-order symmetric tensor $\overline{\grad^s} \gz x$ is cheap to evaluate.
To that end, we cache two scalars that depend on the Jacobian $J$.

\begin{algorithm}[!ht]
  \SetKwInOut{Input}{Given}
  \SetKwInOut{Output}{Return}
  \Input{Source FE vector $\gz x$,
  cached $\gz \tau/J$ and $\boldsymbol{\mathcal{C}}$ for each cell and quadrature point,
  evaluated based on $\overline{\gz u^h}$}
  \Output{action of the FE tangent operator \eqref{eq:algebraic_tangent} on $\gz x$}
  \ForEach{ element $K \in \Omega^h$ }{
        extract local vector values on this element $\gz x_K$ \;
        evaluate the following tensors at each quadrature point using sum factorization :\\
        \enskip$\gz g_s := \overline{\grad^s} \gz x_K$ \tcp*{2nd order symmetric}
        \enskip$\gz g := \overline{\grad} \gz x_K$ \tcp*{2nd order}
        \ForEach{quadrature point $q$ on $K$}{
          evaluate $\boldsymbol{\mathcal{G}} \gz g := \gz g \cdot \left[\gz \tau/J\right]$ \tcp*{2nd order}
          queue $\boldsymbol{\mathcal{G}} \gz g$ for contraction $\overline{\grad} \gz N_i : \boldsymbol{\mathcal{G}} \gz g$\;
          evaluate $\boldsymbol{\mathcal{C}}\gz g_s := \boldsymbol{\mathcal{C}}\gz : \gz g_s$ \tcp*{2nd order symmetric}
          queue $\boldsymbol{\mathcal{C}}\gz g_s$ for contraction $\overline{\grad^s} \gz N_i : \boldsymbol{\mathcal{C}}\gz g_s$ \;
        }
        evaluate queued contractions using sum factorization \;
        distribute results to the destination vector
  }
  \caption{Matrix-free tangent operator: cache material part of the fourth-order spatial tangent stiffness tensor $\boldsymbol{\mathcal{C}}$ and Kirchhoff stress $\gz \tau$.}
  \label{alg:mf_tensor4}
\end{algorithm}

Finally, Algorithm \ref{alg:mf_tensor4} represents the most general case which does not assume any form of the material part of the fourth-order spatial tangent stiffness tensor.
In this case, we have to cache the fourth-order symmetric tensor $\boldsymbol{\mathcal{C}}$ and the Kirchhoff stress $\gz \tau$ for each element and quadrature point, and perform the double contraction with the second-order symmetric tensor $\overline{\grad^s} \gz x$ on-the-fly.

Note that the SIMD vectorization
in \texttt{deal.II} \cite{dealII90}
is applied at the finite element level, i.e. the matrix-free operator is applied simultaneously on several elements (called ``blocks'').
The main reason is that, apart from the point-wise computation of the stresses and elastic tangents, the operations are typically the same on all elements\footnote{%
In general, not all quadrature points in each element group are endowed with identical, non-dissipative constitutive laws that follow the same code path during evaluation or the kinetic quantities and tangents.
This motivates the caching strategy employed in Algorithm 3, as the pre-caching of the chosen data ensures a context in which SIMD parallelization of subsequent operations is guaranteed.
Algorithms 1 and 2 sacrifice some of this generality for decreased memory requirements.
}.
For a given quadrature point, we simultaneously perform tensor evaluations and contractions on all elements\footnote{%
For the heterogeneous materials considered here,
the elements are first grouped according to the material type (matrix or inclusion) and then the SIMD vectorization is applied for each group of elements.
This functionality is provided by the \texttt{deal.II} library in the \textit{MatrixFree} class.
}.
For example, the deformation gradient $\gz F = \gz I + \Grad \, \overline{\gz u^h}$ (second-order tensor) is evaluated at the specific quadrature point of all elements within the ``block''.
This is achieved using templated number types within the \textit{Tensor} class of the \texttt{deal.II} library.
In particular, \texttt{deal.II} provides the \textit{VectorizedArray} templated generic class that defines a unified interface to a vectorized data type
and overloads basic arithmetic operations based on compiler intrinsics.
The core of the sum factorization algorithms is provided by the
\textit{FEEvaluation} class, that already supports all the necessary contraction (such as $\overline{\grad^s} \gz N_i : \boldsymbol{\mathcal{C}}\gz g_s$) to implement the three algorithms proposed above.
Thanks to the object-oriented design of the \texttt{deal.II} library, the core of those matrix-free algorithms can be implemented in only a handful of lines.
Distributed memory MPI parallelization of the matrix-free operators is done using the standard domain decomposition technique, where each MPI process is responsible for applying the operator only on the so-called locally owned elements.
For more information about the implementation details and data structures for sum factorization approaches in \texttt{deal.II} (including hanging constraints, storage of the inverse Jacobian of the transformation from unit to real cell, MPI and distributed memory parallelization as well SIMD vectorization) and performance tests, see \cite{kronbichler12}.

\section{Geometric multigrid preconditioning}
\label{sec:gmg}

Efficient and scalable matrix-vector products are not enough to obtain a method suitable for large computations. In particular, linear scaling with respect to the number $N$ of unknown DoFs is desirable.
The efficiency of standard iterative solvers and preconditioners deteriorates for large systems of linear equations due to increasing iteration counts under mesh refinement.
In contrast to other preconditioners, geometric multigrid (GMG) preconditioners \cite{Bramble1990, Briggs2000, Janssen2011,May2015} can be proven to result in iteration counts that are independent of the number of mesh refinements by smoothing the residual on a hierarchy of meshes. For further information about geometric multigrid methods we refer the reader to \cite{Briggs2000,Hackbusch1985,Wesseling1992}.
Implementational details of the MPI-parallel multigrid method in \texttt{deal.II} are discussed in \cite{Clevenger2018}.

For this work, we adopt the matrix-free level operators within the GMG.
In general, given a heterogeneous non-linear material, the definition of transfer operators\footnote{Finite dimensional linear operators from a vector space associated with the fine-scale mesh to the coarse-scale mesh, and in the opposite direction.} (for restriction and prolongation) and level operators is not trivial. One approach is to start from an arbitrary heterogeneous material at the fine scale and employ homogenization theory \cite{Suquet1987, Hill1972,Hashin1983,Castaneda1997} to design suitable transfer operators \cite{Miehe2007}\footnote{
  We note that this approach is limited to the case where for each quadrature point within an element the fourth-order tangent stiffness tensor is the same.
}.
In this approach, the fine-scale displacement is decomposed into long-wave and short-wave contributions by an additive split, where the former is associated to the homogeneous response of a patch of elements and the latter represents fluctuations due to the heterogeneous structure of the patch.
However, for a patch of elements consisting of the same material, the transfer operators reduce to the standard geometric transfer.
Based on this observation we adopt the following approach: It is assumed that the mesh at the coarsest level can accurately describe the heterogeneous finite-strain elastic material.
In this case, the restriction and prolongation operations for patches of elements are always performed for the same material and thus admits usage of standard geometric transfer operators.
We define the level operators as follows: Recall that the fine scale tangent operator \eqref{eq:algebraic_tangent} is obtained by linearization around the current displacement field $\gz u^h$. We then restrict this displacement field to all multigrid levels $\{l\}$ and evaluate tangent operators using the matrix-free algorithms from Section \ref{sec:mf}.
Essentially, level operators correspond to the linearization around the smoothed representation of displacement field\footnote{
Note that this differs from the classical approach where the (tangent) operator $\gz A^{l+1}$ on level $l+1$ is directly related to the (tangent) operator $\gz A^{l}$ on level $l$ via $\gz A^{l+1}=\gz I^{l+1}_{l} \gz A^l \gz I^l_{l+1}$, where $\gz I^l_{l+1}$ and $\gz I^{l+1}_l$ are global prolongation (coarse-to-fine) and restriction (fine-to-coarse) operators, respectively.
}.

To summarize, in this section we propose a method by which to construct level tangent operators for finite-strain elasticity to be employed within the GMG preconditioner.
The other parts of the GMG preconditioner such as level smoothers, matrix-free interlevel transfer operators and the coarse level iterative solver are provided by the \texttt{deal.II} library.
The efficiency of the preconditioner with the considered level operators will be assessed by numerical examples in the next section.
We will, however, refrain from studying the machine performance aspects of the multigrid preconditioner as a whole as the focus of this paper is the matrix-free implementation of tangent operators of finite-strain elasticity both at the fine mesh level as well as the coarser multigrid levels.
Finally, we note that at the time of writing, a $p-$ version of the GMG was not available within the \texttt{deal.II} library and therefore not considered here.

For examples of other approaches that combine homogenization and multigrid methods, see \cite{Bayreuther2003,Fish1995} for applications using small strain elastic materials and \cite{Kaczmarczyk2010} with applications in fracture.

\section{Numerical examples}
\label{sec:example}

\begin{figure}[!ht]
  \centering
  \begin{subfigure}[b]{0.45\textwidth}
      \centering
      \includegraphics[width=\textwidth]{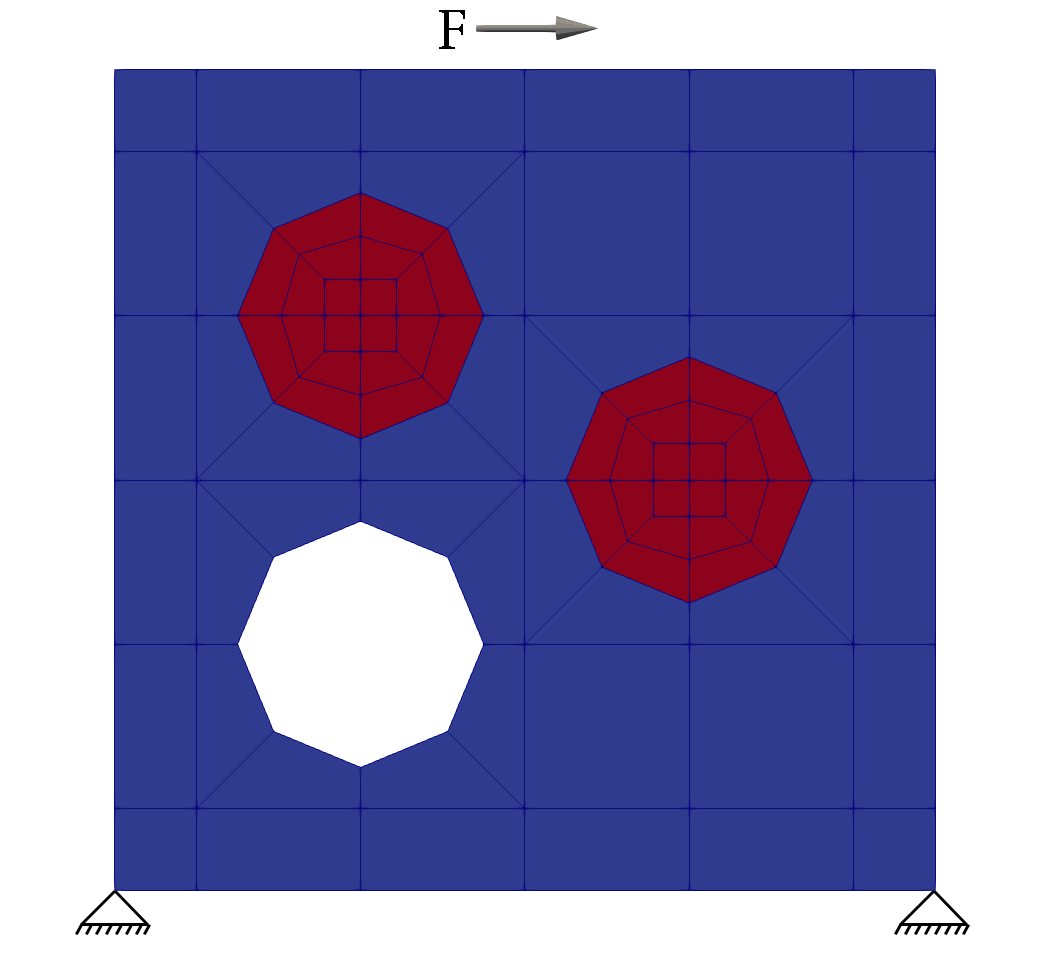}
      \caption{2D coarse mesh}
  \end{subfigure}
  \begin{subfigure}[b]{0.45\textwidth}
    \centering
    \includegraphics[width=\textwidth]{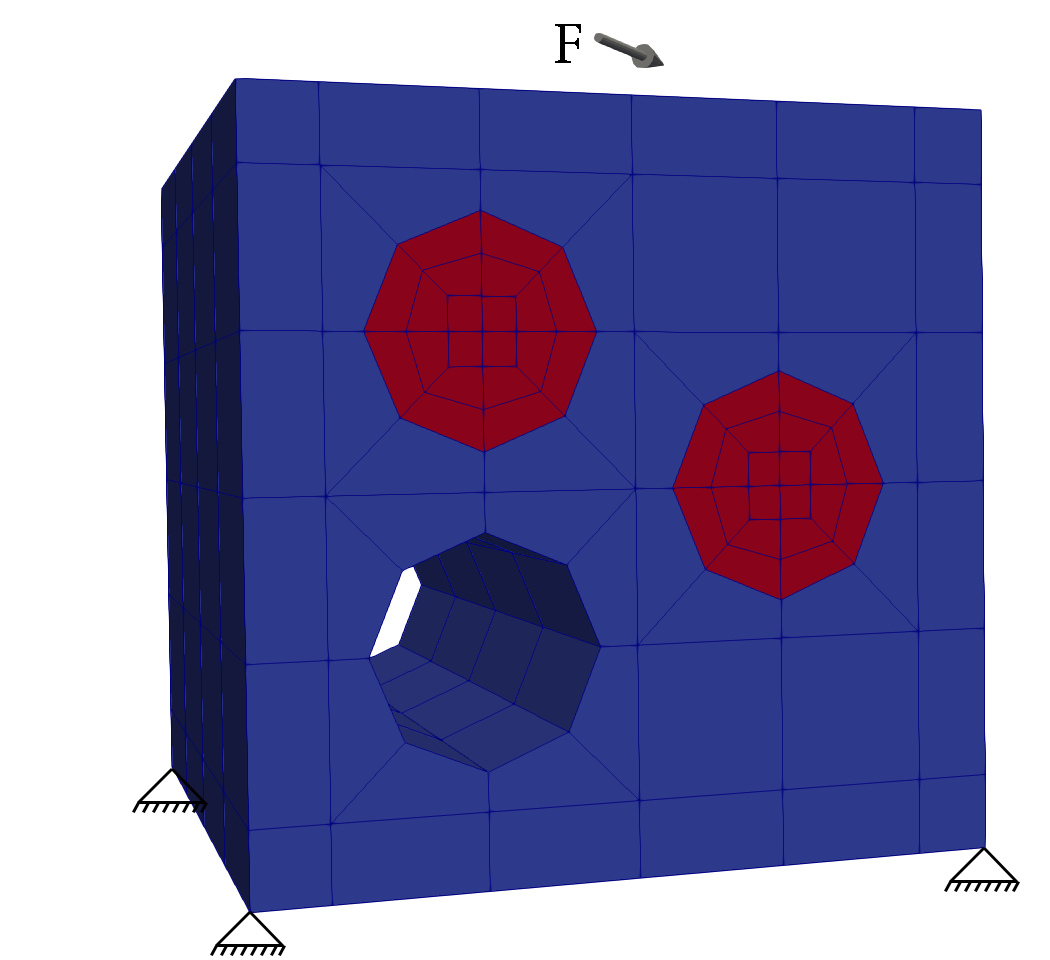}
    \caption{3D coarse mesh}
  \end{subfigure}
  \caption{Discretization of the heterogeneous material at the coarsest mesh level and the prescribed boundary conditions.}%
  \label{fig:miehe}
\end{figure}

\begin{figure}[!ht]
  \centering
  \begin{subfigure}[b]{0.45\textwidth}
      \centering
      \includegraphics[width=\textwidth]{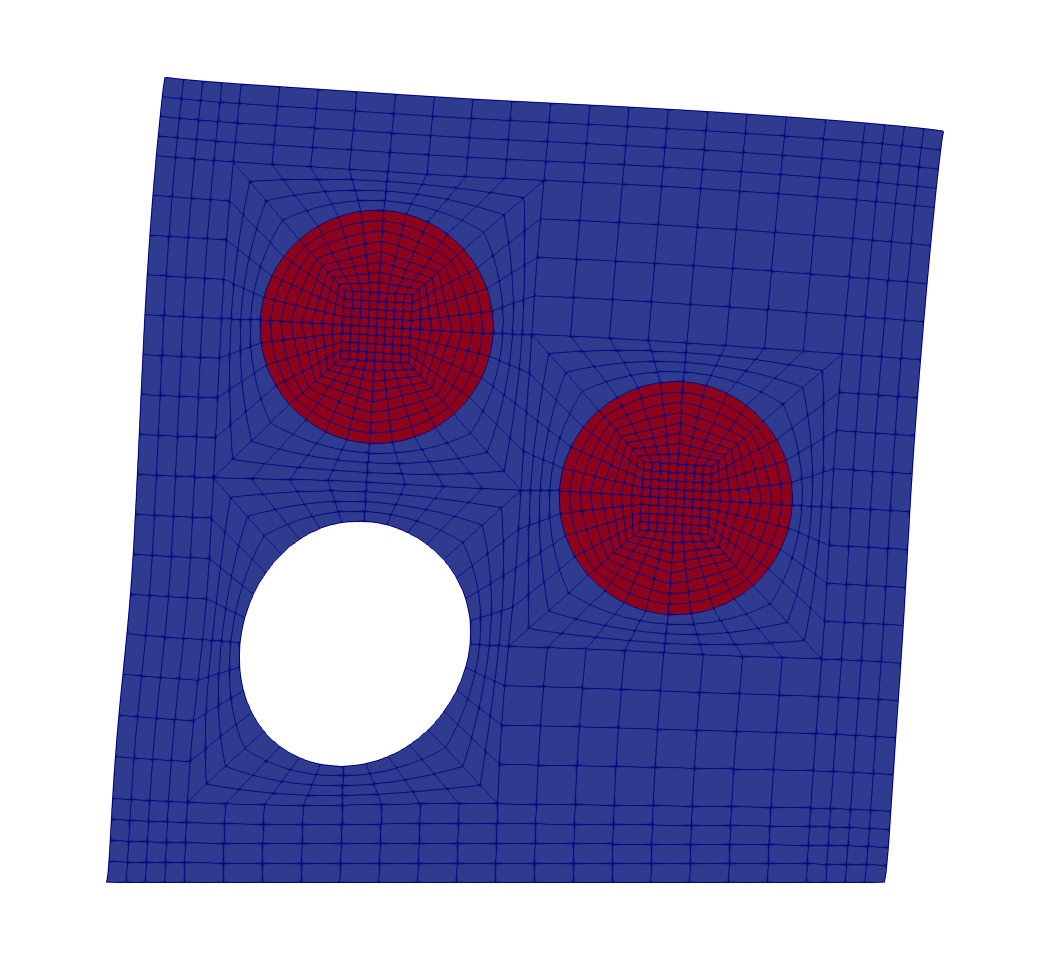}
      \caption{2D deformed mesh}
  \end{subfigure}
  \begin{subfigure}[b]{0.45\textwidth}
    \centering
    \includegraphics[width=\textwidth]{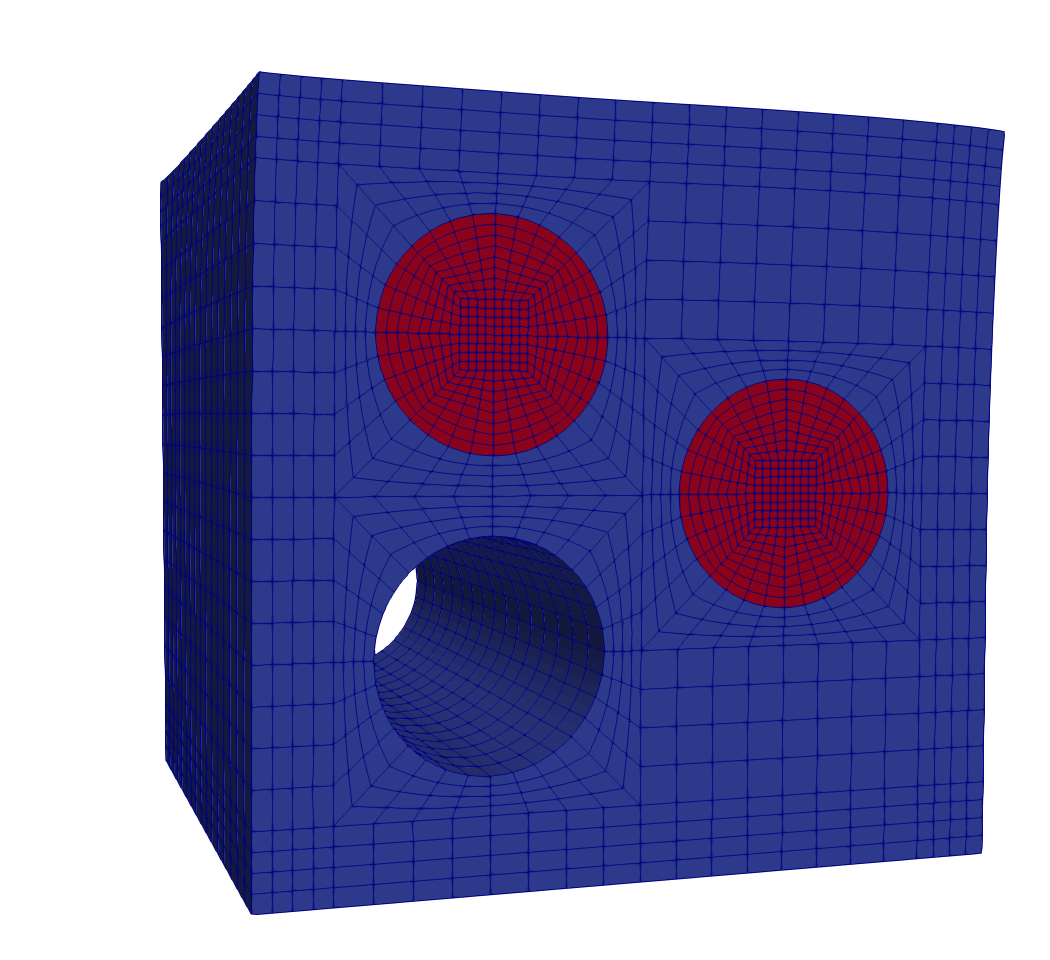}
    \caption{3D deformed mesh}
  \end{subfigure}
  \caption{Deformed meshes at the final loading step with quadratic elements and two global mesh refinements.}%
  \label{fig:miehe_deformed}
\end{figure}

In this section, we apply the proposed matrix-free operator algorithms for finite-strain hyperelastic materials as well as the geometric multigrid preconditioner to a benchmark problem from \cite{Miehe2007}. The coarse meshes for the 2D and 3D problems are illustrated in Figure \ref{fig:miehe}. The 2D material consists of a square (depicted in blue), a hole and two spherical inclusions (depicted in red). The size of the square domain is $10^{-3} \rm{mm}$.
The matrix material is taken to have Poisson's ratio $0.3$ and shear modulus $\mu=0.4225 \times 10^6\, \rm{N/mm^2}$. The inclusion is taken to be $100$ times stiffer.
The 3D material is obtained by extrusion of the 2D geometry into the third dimension.
Note that although a relatively coarse mesh is used at the coarsest multigrid level, the geometry of the inclusions is captured accurately by manifold descriptions of the boundaries and interfaces.
Both domains are fully fixed along the bottom surface and a distributed load is applied at the top in the $(1,0)$ or $(1,1,0)$ direction for the 2D and the 3D problems, respectively.
The force density $12.5 \times 10^3 \rm{N/mm^2}$ or $12.5 \sqrt{2} \times 10^3 \rm{N/mm^3}$ is applied in 5 steps for the 2D and the 3D problems, respectively.
With respect to the nonlinear solver, the displacement tolerance of the $\mathcal{\ell}_2$ norm for the Newton update is taken to be $10^{-5}$ whereas the relative and absolute tolerances for the residual forces is $10^{-8}$. The relative convergence criteria for the linear solver is $10^{-6}$.
Figure \ref{fig:miehe_deformed} shows deformed meshes at the final loading step with quadratic elements and two global mesh refinements.
\begin{figure}[!ht]
  \centering
  \begin{subfigure}[b]{0.8\textwidth}
    \centering
    \includegraphics[width=\textwidth]{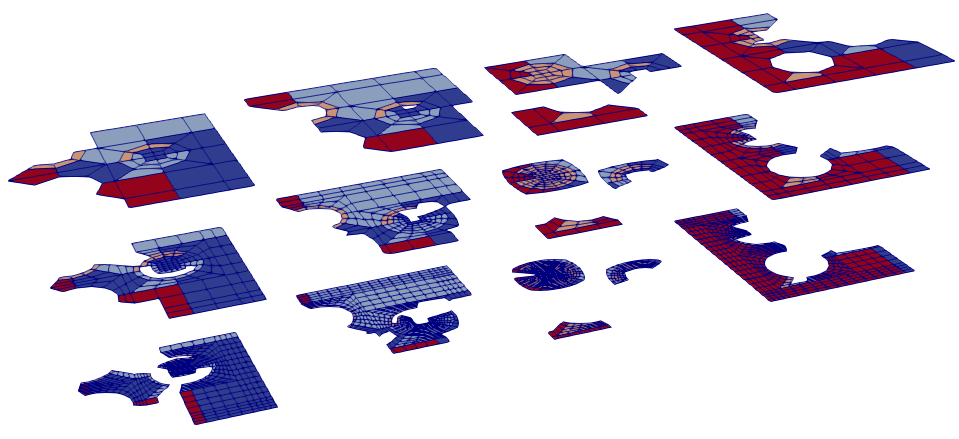}
  \end{subfigure}
  \caption{2D multigrid mesh after two global refinements distributed into four MPI processes. The color indicates the process owning the element.}%
  \label{fig:miehe_gmg}
\end{figure}
The same distributed approach to mesh decomposition (as most evident on the finest level) is adopted on the multigrid levels.
Figure \ref{fig:miehe_gmg} illustrates the hierarchy of multigrid meshes for the 2D problem after two global refinements for the mesh introduced in figure \ref{fig:miehe_fine_level}.
\begin{table}
  \centering
  \begin{tabular}{ccccc}
  \hline
    $p$ & $q$ & $N_{gref}$ & $N_{el}$ & $N_{DoF}$ \\
  \hline
    1 & 2 & 7 & 1441792 & 2887680 \\
    2 & 3 & 6 & 360448 & 2887680 \\
    3 & 4 & 5 & 90112 & 1625088 \\
    4 & 5 & 5 & 90112 & 2887680 \\
    5 & 6 & 5 & 90112 & 4510720 \\
    6 & 7 & 4 & 22528 & 1625088 \\
    7 & 8 & 4 & 22528 & 2211328 \\
    8 & 9 & 4 & 22528 & 2887680 \\
  \hline
  \end{tabular}
  \caption{Parameters for the 2D benchmark: $p$ is the polynomial degree,
  $q$ is the number of quadrature points in 1D, $N_{gref}$ is the number of global mesh refinements, $N_{el}$ is the number of elements and $N_{DoF}$ is the number of DoFs.
  }
  \label{tab:input_parameters_2d}
\end{table}

\begin{table}
  \centering
  \begin{tabular}{ccccc}
  \hline
    $p$ & $q$ & $N_{gref}$ & $N_{el}$ & $N_{DoF}$ \\
  \hline
    1 & 2 & 4 & 1441792 & 4442880 \\
    2 & 3 & 3 & 180224 & 4442880 \\
    3 & 4 & 2 & 22528 & 1891008 \\
    4 & 5 & 2 & 22528 & 4442880 \\
  \hline
  \end{tabular}
  \caption{Parameters for the 3D benchmark: $p$ is the polynomial degree,
  $q$ is the number of quadrature points in 1D, $N_{gref}$ is the number of global mesh refinements, $N_{el}$ is the number of elements and $N_{DoF}$ is the number of DoFs.
  }
  \label{tab:input_parameters_3d}
\end{table}

Aligned with the approach taken in similar previous investigations, compare e.g. \cite{kronbichler2017fast},
we limit ourselves to comparisons made on a node level in this study.
Examinations performed in this manner should already indicate the differences between the matrix-free and matrix-based approaches, and give some preliminary understanding of the scalability potential of the numerical implementation without confounding factors related to data exchange through interconnects via MPI.
The benchmark calculations are performed for 2D and 3D problems with various combinations of polynomial degrees and number of global refinements, as is stated in Tables \ref{tab:input_parameters_2d} and \ref{tab:input_parameters_3d}.
We chose a problem size around $10^6$ DoFs for each test case as this was an upper limit due to the memory requirements of the matrix-based approach with a high polynomial order finite element.
Even so, the problem size is large enough so that the sparse matrix (in excess of $1$ Gb when using the lowest polynomial order finite element discretization) will not fit into the L3 cache (on the order of tens of Mb in size), and therefore we can observe that the matrix-based approach is memory bound.

We conduct two performance studies using the described discretization, the first being for only matrix-vector multiplication, and the second for the preconditioned iterative linear solver.
The results of the matrix-free approach with a geometric multigrid preconditioner are compared to the matrix-based approach in conjunction with the algebraic multigrid (AMG) preconditioner using the package ML \cite{Gee2006a} of the Trilinos \cite{Heroux2005} library version 12.12.1.
The aggregation threshold for the AMG preconditioner is taken as $10^{-4}$.
The geometric multigrid algorithm uses a V-cycle with two pre- and post-smoothing steps using Chebyshev polynomials \cite{Varga2009} of fourth degree, provided
by the \texttt{deal.II} library via the \textit{PreconditionChebyshev} class.
The largest eigenvalue $\lambda_{\rm{max}}$ of the level operators is estimated from $30$ iterations of the CG solver and the smoother
is set to target the range $[0.6\lambda_{\rm{max}}, 1.2\lambda_{\rm{max}}]$.
On the coarsest level, we use the Chebyshev preconditioner as a solver \cite{Varga2009} to reduce the residual by three orders of magnitude.
Since this smoother only uses diagonal elements of the operator it is readily compatible with the matrix-free approach while most default smoothers rely on an explicit matrix form.

Computations were performed on a single node of two clusters:
A node on the ``Emmy'' cluster at RRZE, FAU has two Xeon 2660v2 ``Ivy Bridge'' chips (10 cores per chip + SMT) running at 2.2 GHz with 25 MB shared cache per chip and 64 GB of RAM. Intel's compiler version 18.03 with flags ``-O3 -march=native'' and Intel MPI version 18.03 are used.
The ``IWR'' results were obtained on a single machine with eight Xeon E7-8870 ``Sandy Bridge'' chips (10 cores per chip  + SMT) running at 2.4 GHz with 30 MB shared cache per chip. For the simulation, GCC version 8.1.0 with flags ``-O3 -march=native'' and OpenMPI version 3.0.0 were used.
Unless noted otherwise, the simulations are performed with 20 MPI processes, i.e. one fully occupied node.

\subsection{Matrix-vector multiplication}

\begin{figure}[!ht]
  \centering
  \begin{subfigure}[b]{0.49\textwidth}
      \centering
      \includegraphics[width=\textwidth]{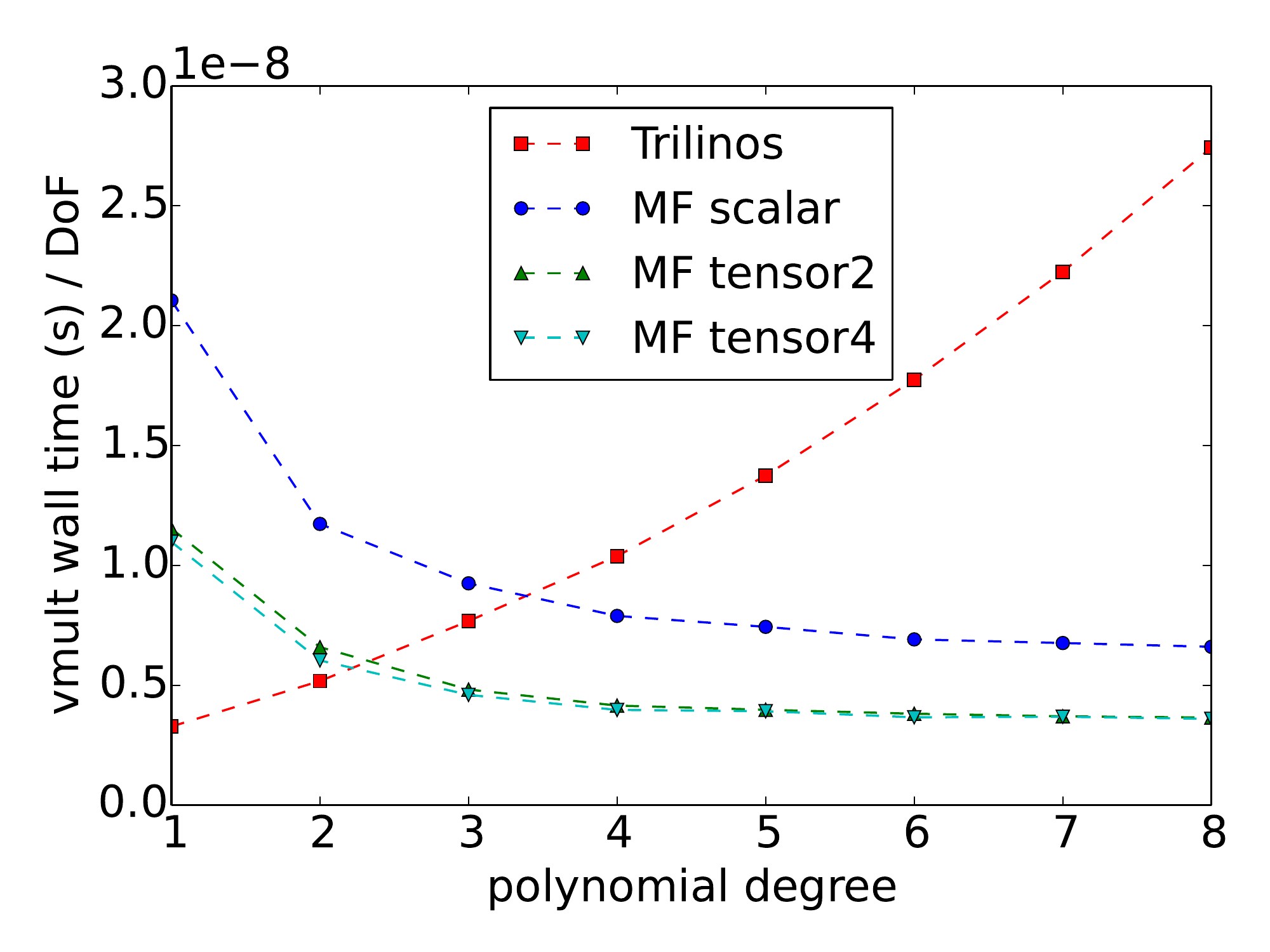}
      \caption{matrix-vector product (2D)}
      \label{fig:benchmark_miehe_Emmy_vmult2}
  \end{subfigure}
  \begin{subfigure}[b]{0.49\textwidth}
    \centering
    \includegraphics[width=\textwidth]{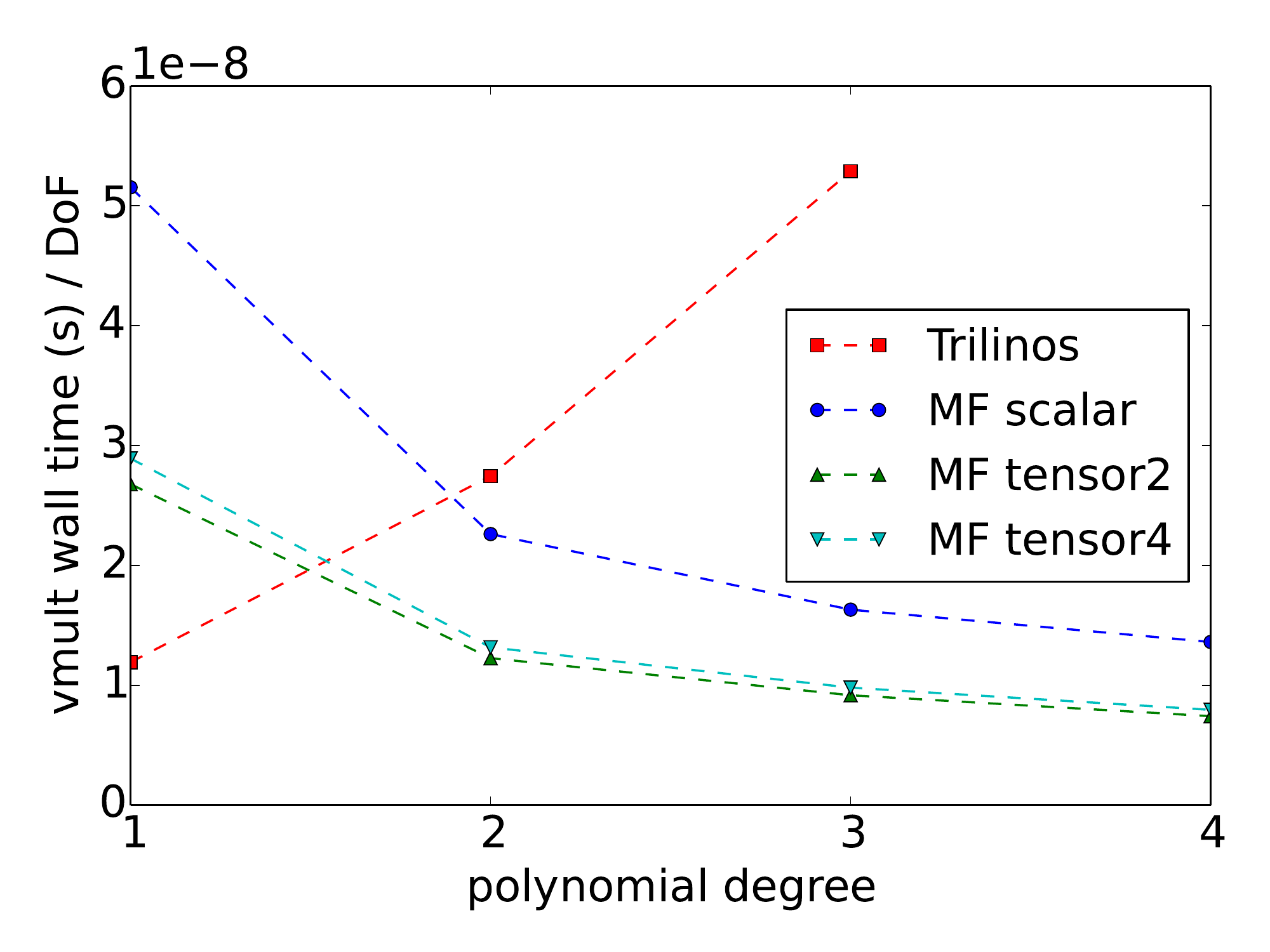}
    \caption{matrix-vector product (3D)}
    \label{fig:benchmark_miehe_Emmy_vmult3}
  \end{subfigure}
  ~
  \begin{subfigure}[b]{0.49\textwidth}
      \centering
      \includegraphics[width=\textwidth]{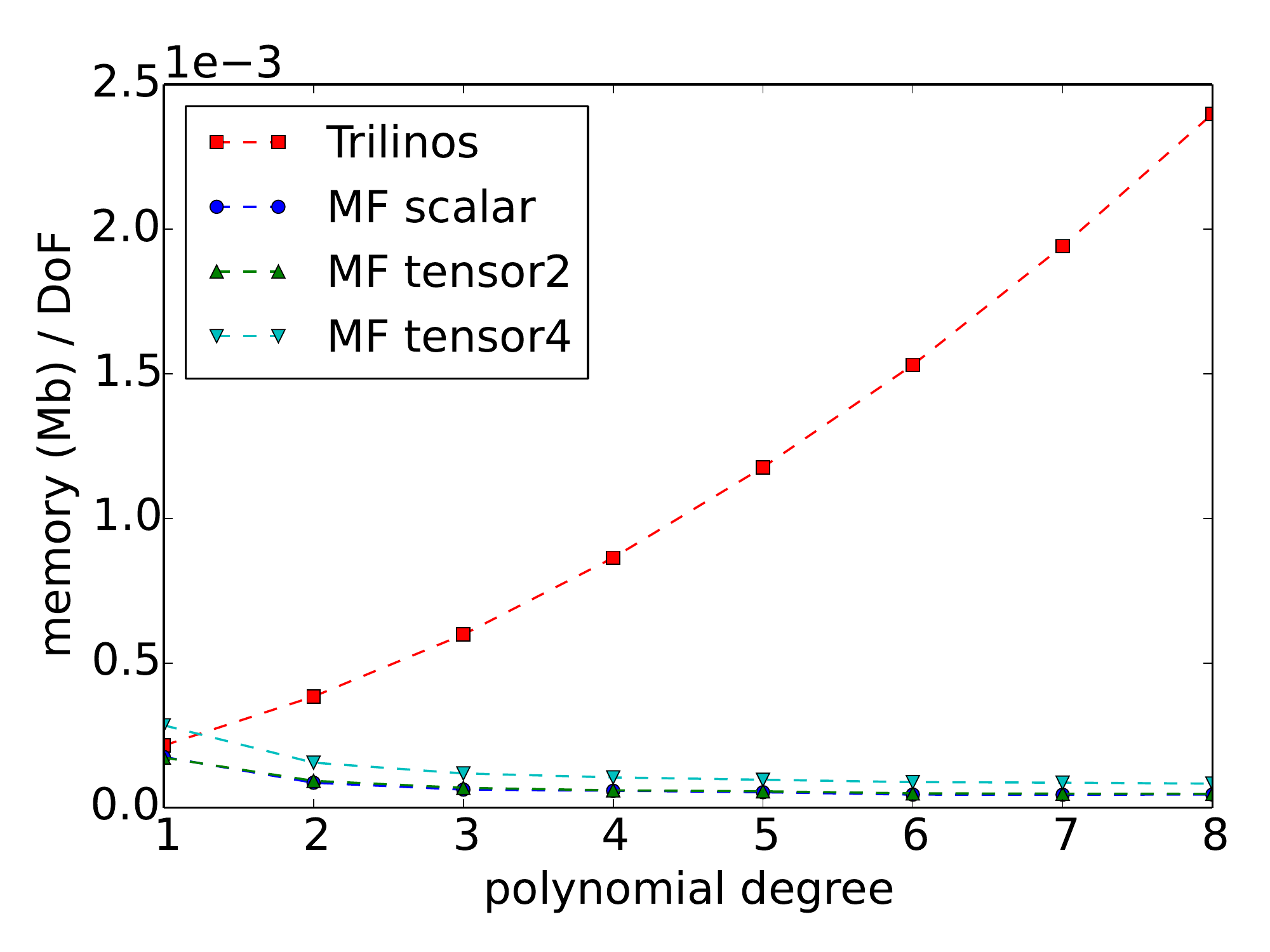}
      \caption{memory consumption (2D)}
      \label{fig:benchmark_miehe_Emmy_memory2}
  \end{subfigure}
  \begin{subfigure}[b]{0.49\textwidth}
    \centering
    \includegraphics[width=\textwidth]{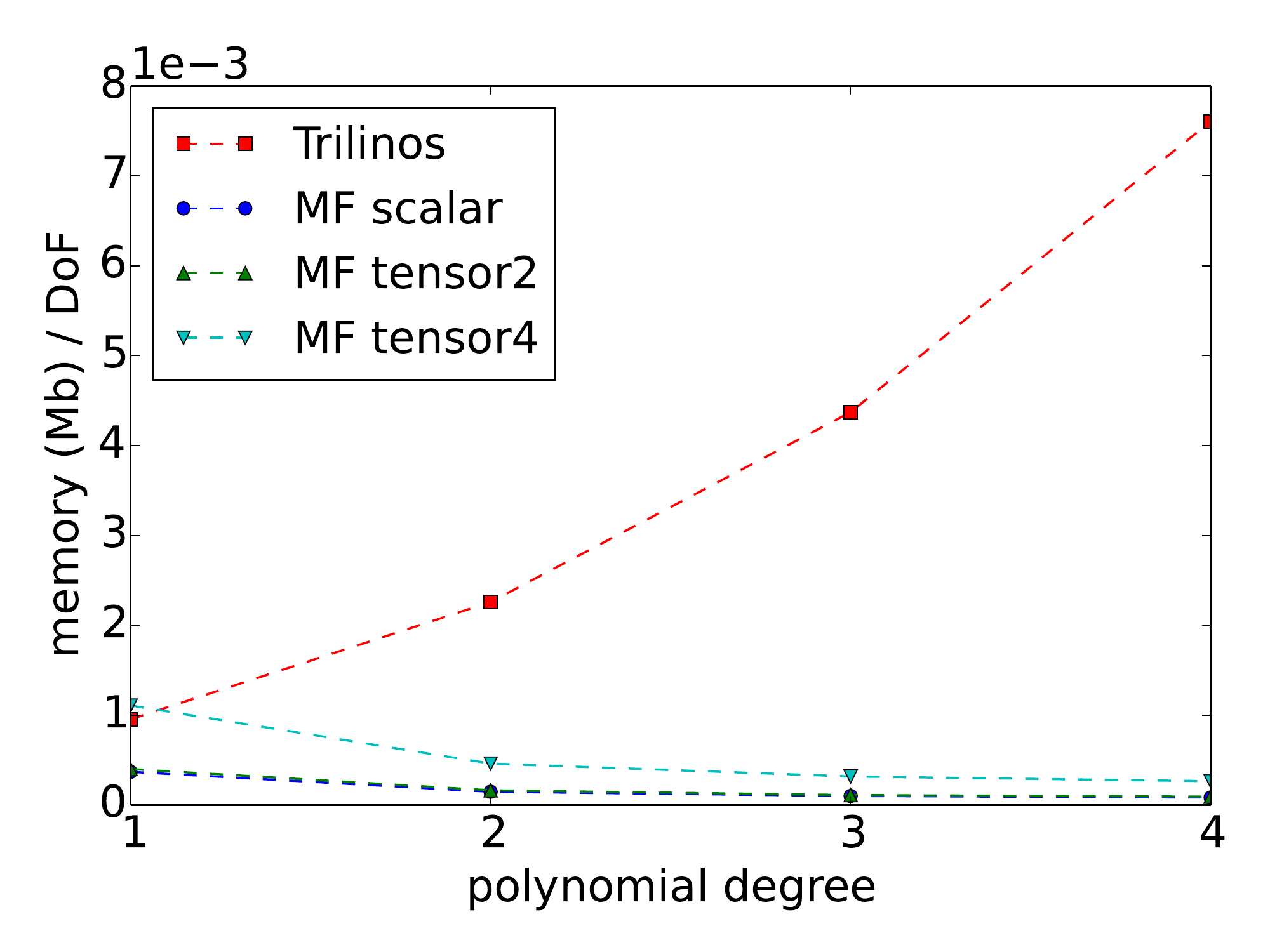}
    \caption{memory consumption (3D)}
    \label{fig:benchmark_miehe_Emmy_memory3}
  \end{subfigure}
  \caption{Emmy cluster, RRZE, matrix-vector multiplication.}%
  \label{fig:benchmark_miehe_Emmy}
\end{figure}

\begin{figure}[!ht]
  \centering
  \begin{subfigure}[b]{0.49\textwidth}
      \centering
      \includegraphics[width=\textwidth]{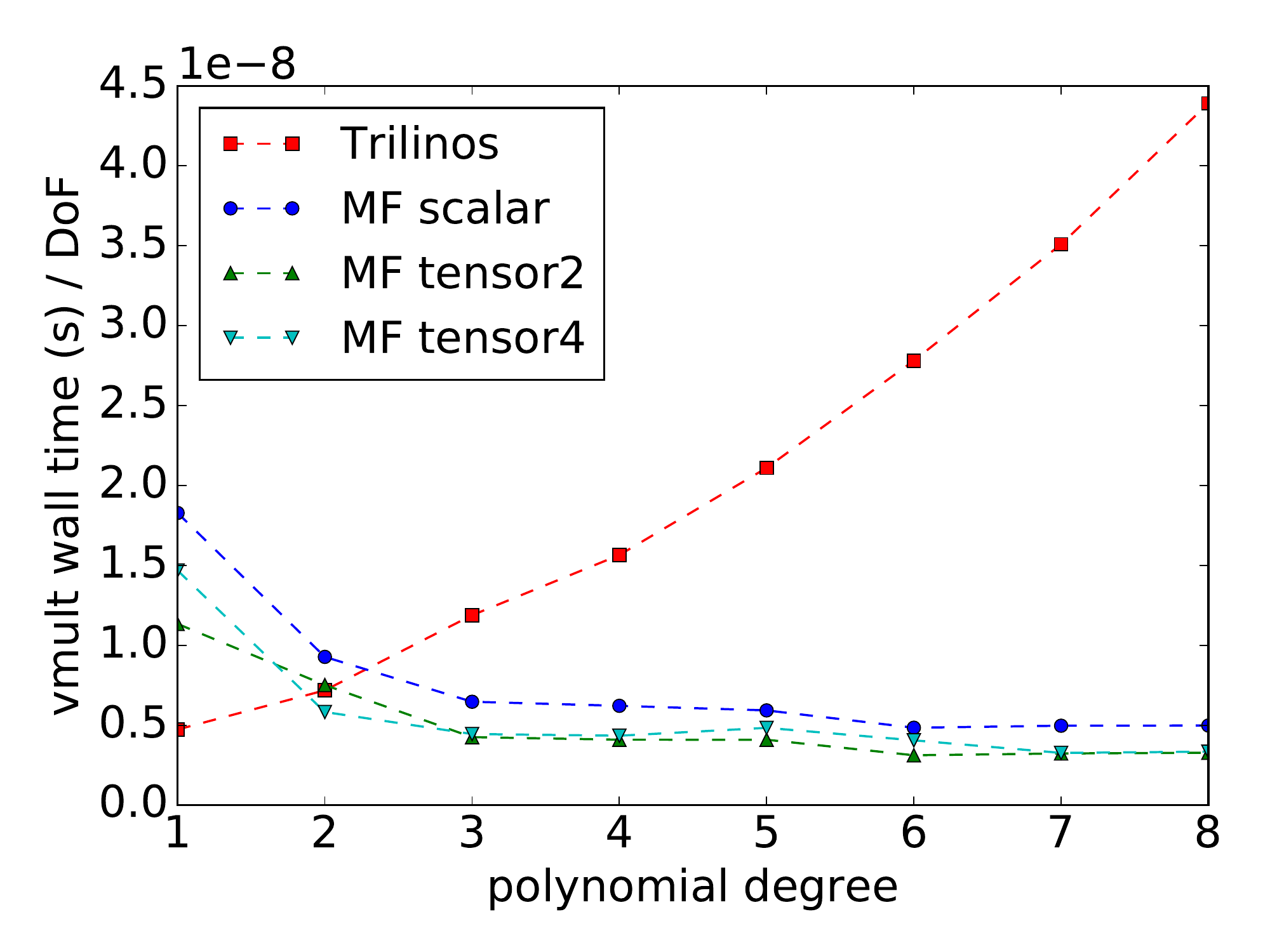}
      \caption{matrix-vector product (2D)}
      \label{fig:benchmark_miehe_IWR_vmult2}
  \end{subfigure}
  \begin{subfigure}[b]{0.49\textwidth}
    \centering
    \includegraphics[width=\textwidth]{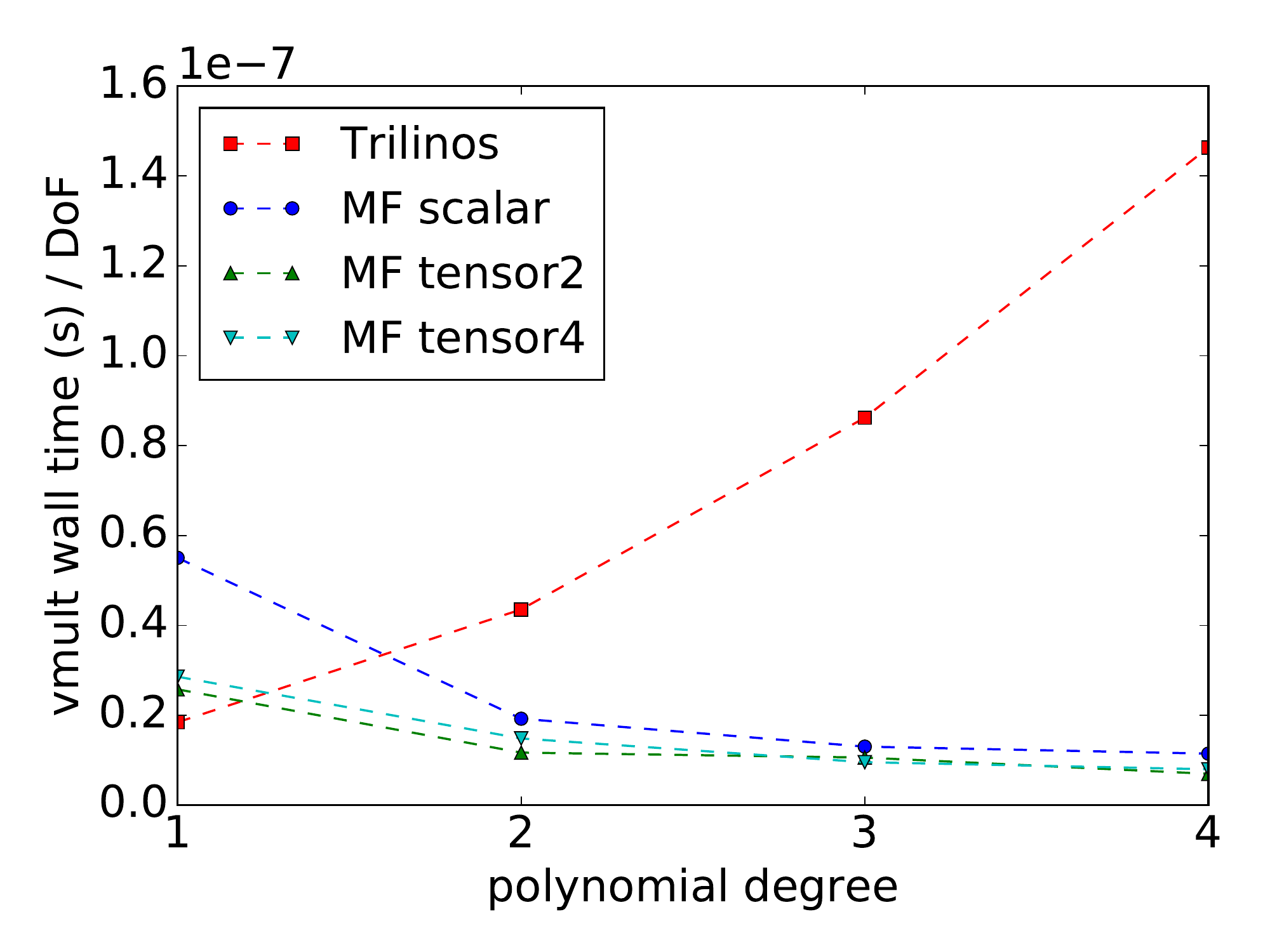}
    \caption{matrix-vector product (3D)}
    \label{fig:benchmark_miehe_IWR_vmult3}
  \end{subfigure}
  ~
  \begin{subfigure}[b]{0.49\textwidth}
      \centering
      \includegraphics[width=\textwidth]{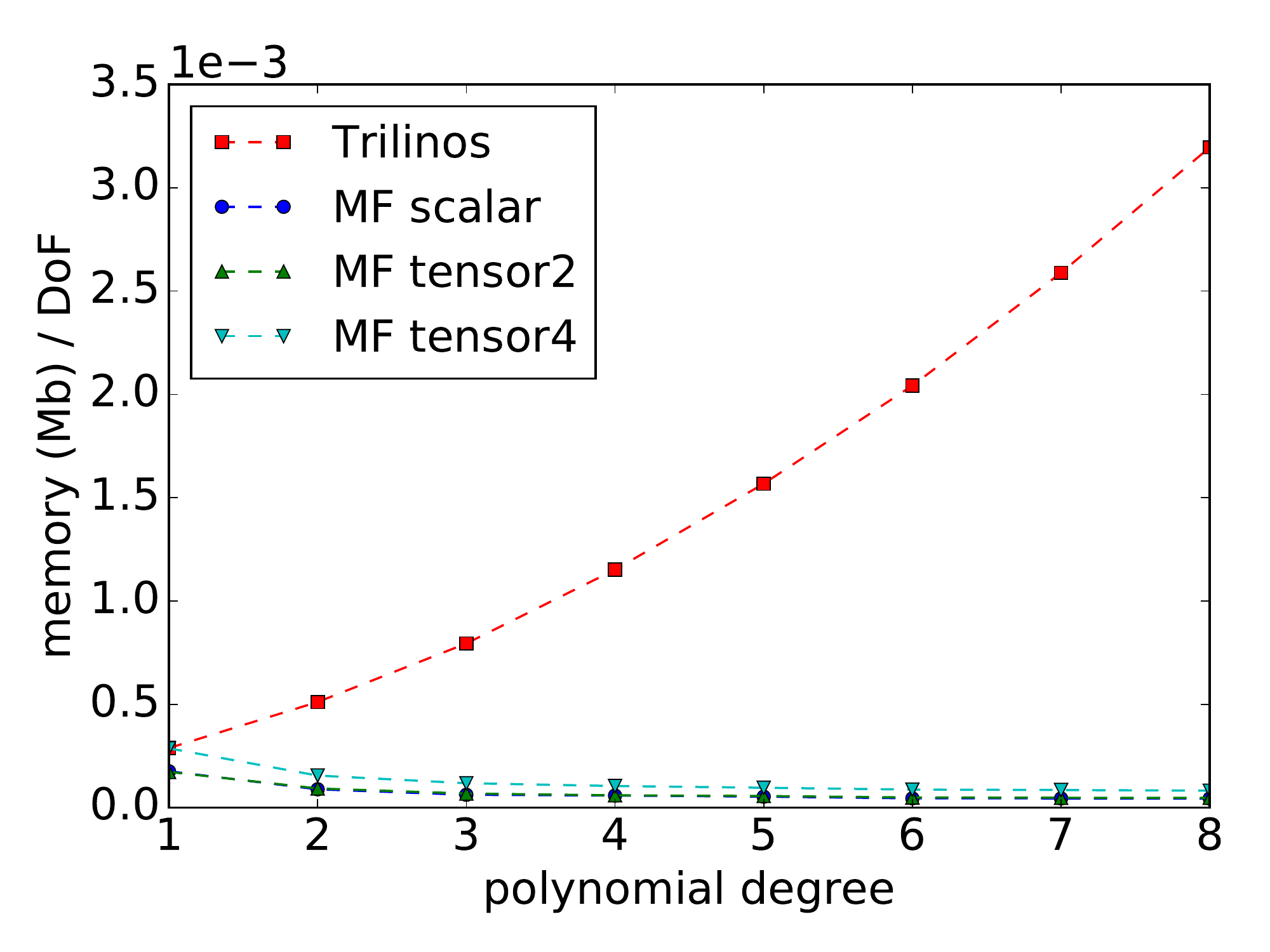}
      \caption{memory consumption (2D)}
      \label{fig:benchmark_miehe_IWR_memory2}
  \end{subfigure}
  \begin{subfigure}[b]{0.49\textwidth}
    \centering
    \includegraphics[width=\textwidth]{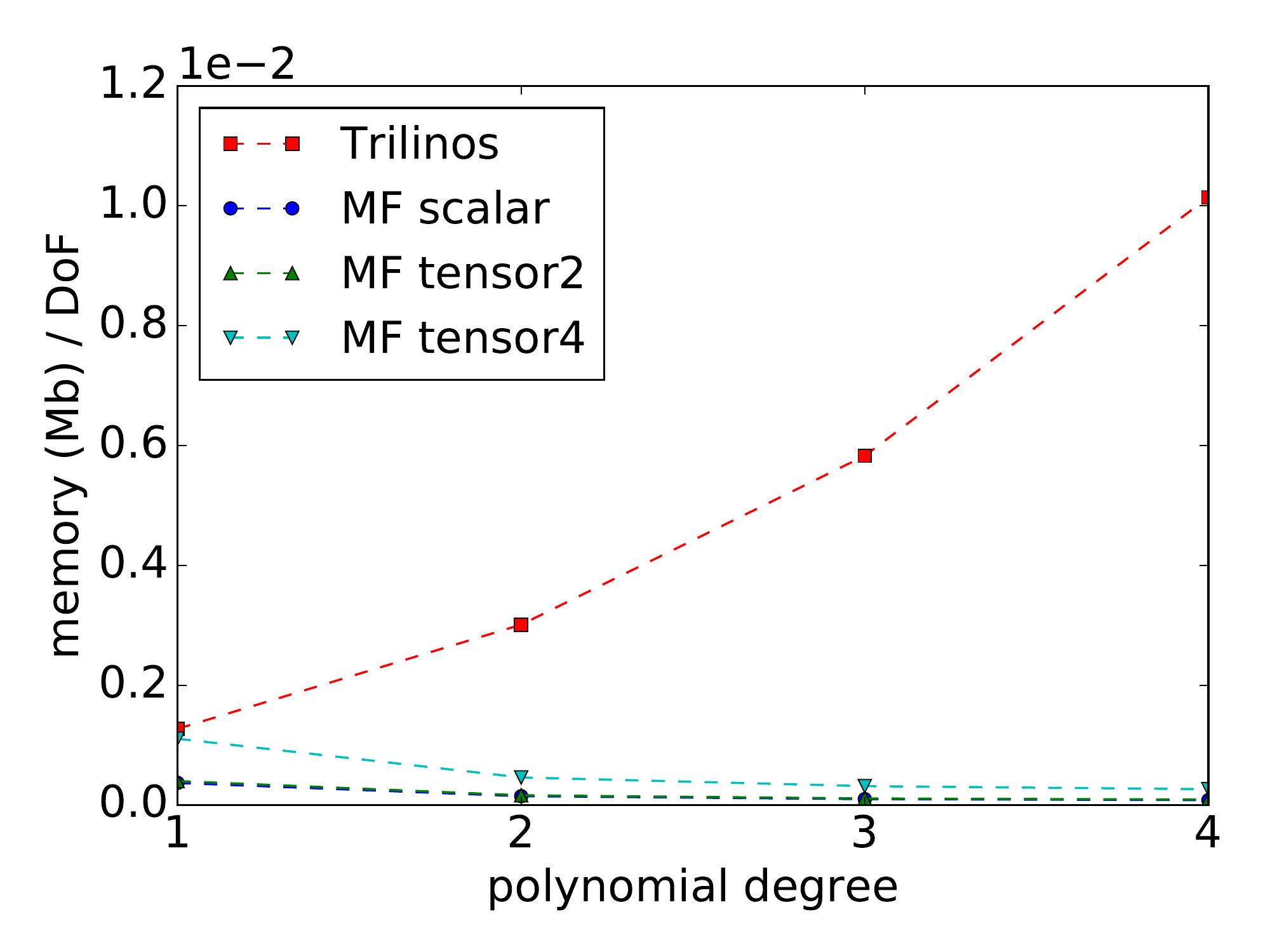}
    \caption{memory consumption (3D)}
    \label{fig:benchmark_miehe_IWR_memory3}
  \end{subfigure}
  \caption{IWR cluster, matrix-vector multiplication}%
  \label{fig:benchmark_miehe_IWR}
\end{figure}

Figures \ref{fig:benchmark_miehe_Emmy} and \ref{fig:benchmark_miehe_IWR} show the results of the matrix-vector multiplication for the considered finite-strain hyperelastic benchmark problem as performed on the two clusters.
We are interested in the wall-clock time (average over 10 runs per each linear solver step) per DoF as a first metric.
As expected, the matrix-vector multiplication becomes very expensive for higher order discretization for sparse matrix-based approaches,
as the matrix becomes more and more dense.
The performance results for finite-strain elasticity operators considered in this study are qualitatively similar to those reported for the Laplace operator in \cite{kronbichler12}.
Note that we consider only non-Cartesian meshes together with additional mappings (namely, manifold descriptions of the boundaries and interfaces \cite{Arndt2017a} to capture the geometry as well as
the linearization-related mapping required for evaluation of gradients with respect to the deformed configuration).
Both factors increase the memory consumption and make element operations more expensive.
Nonetheless,
the matrix-free implementation is faster than the matrix-based counterpart already for bi-cubic Lagrangian basis in 2D and tri-quadratic in 3D.
Figures \ref{fig:benchmark_miehe_Emmy_vmult2}, \ref{fig:benchmark_miehe_Emmy_vmult3}, \ref{fig:benchmark_miehe_IWR_vmult2} and \ref{fig:benchmark_miehe_IWR_vmult3} clearly show the influence of the different caching strategies,
described in Algorithms \ref{alg:mf_scalar}, \ref{alg:mf_tensor2} and \ref{alg:mf_tensor4}.
Two conclusions can be drawn from these results:
The scalar caching implementation is the most time consuming of the three, as at each step we additionally evaluate the gradient of the displacement field in the referential configuration that is required to evaluate the Kirchhoff stress $\gz \tau$.
There is little difference in terms of the wall-clock time per DoF between the approach where the material part of the fourth-order spatial tangent stiffness tensor is cached (``tensor4'') and the one where we only cache the second-order Kirchhoff tensor and utilize the chosen hyperelastic constitutive equation to efficiently implement the action of the material part of the fourth-order spatial tangent stiffness tensor on the second-order symmetric tensor (``tensor2'').
This indicates that the time savings we get from the latter approach are insignificant within the complete matrix-free operator implementation.
Consequently, we can apply the matrix-free approach to any finite-strain material model by caching the material part of the fourth-order spatial tangent stiffness tensor in addition to the Kirchhoff stress and expect this to be faster than the matrix-based implementation.
It is therefore feasible to define highly performant generic operators that are independent of any applied constitutive laws,
as long as the resulting tangent operator is expressed via \eqref{eq:algebraic_tangent}.

As outlined in the introduction, it is not only the performance of the matrix-free vector multiplication, but also the memory requirement to store a sparse matrix which is the driving force behind matrix-free methods.
Figures \ref{fig:benchmark_miehe_Emmy_memory2}, \ref{fig:benchmark_miehe_Emmy_memory3}, \ref{fig:benchmark_miehe_IWR_memory2} and \ref{fig:benchmark_miehe_IWR_memory3} demonstrate that even with caching one fourth-order symmetric tensor and one second-order tensor at each quadrature point (``tensor4''), the matrix-free approach takes much less memory than its matrix-based counterpart.
In those graphs, we report the overall memory needed to cache additional quantities required for the tangent operator here considered,
as well as the memory required by the \texttt{deal.II} \textit{MatrixFree} object, which caches the inverse of the Jacobian of the transformation
from unit to real cell \cite{kronbichler12}.
Note that for Algorithm \ref{alg:mf_scalar} we need to evaluate gradients with respect to both deformed and undeformed
configurations, and therefore two inverses will be cached via two different \textit{MatrixFree} objects.

Finally, we analyze the node-level performance of the two approaches using the roofline model.
Memory bandwidth (MBytes/s) and the number of floating point operations per second (MFLOP/s) are measured using the \texttt{MEM\_DP} group of the LIKWID \cite{likwid} tool, version 4.2.1.
The measurements are done on a single socket of RRZE cluster by pinning all 10 MPI processes to this socket, whose
memory bandwidth\footnote{As measured with the \textit{load\_avx} benchmark of \textit{likwid-bench} \cite{likwid}.} and peak performance are $B=47.2 \,\rm{GB/s}$ and $P=176\, \rm{GFlop/s}$, respectively.
Turbo mode of the CPU was disabled.

\begin{figure}[!ht]
  \centering
  \begin{subfigure}[b]{0.49\textwidth}
    \centering
    \includegraphics[width=\textwidth]{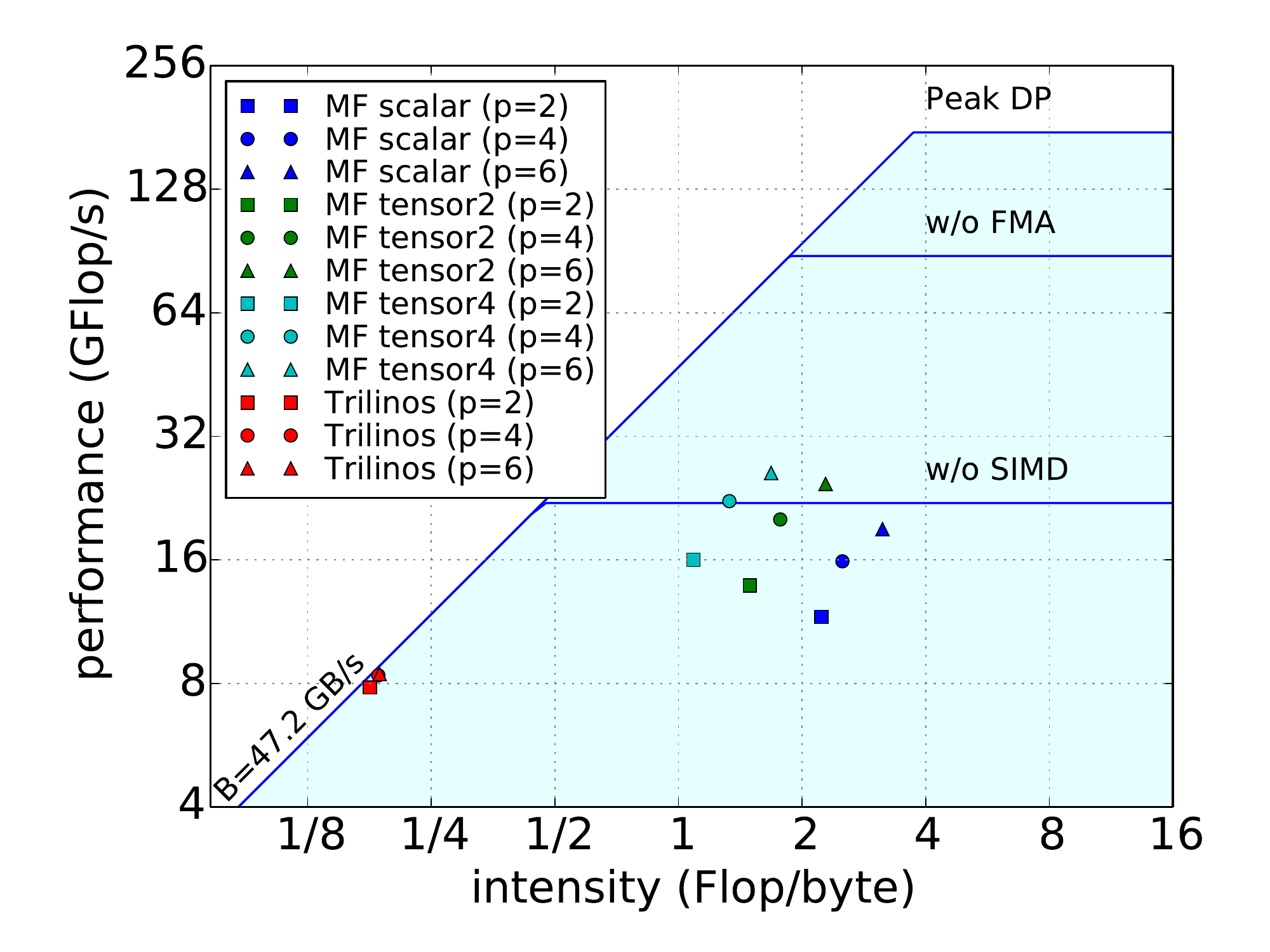}
    \caption{2D}
    \label{fig:roofline_2d}
  \end{subfigure}
  \begin{subfigure}[b]{0.49\textwidth}
    \centering
    \includegraphics[width=\textwidth]{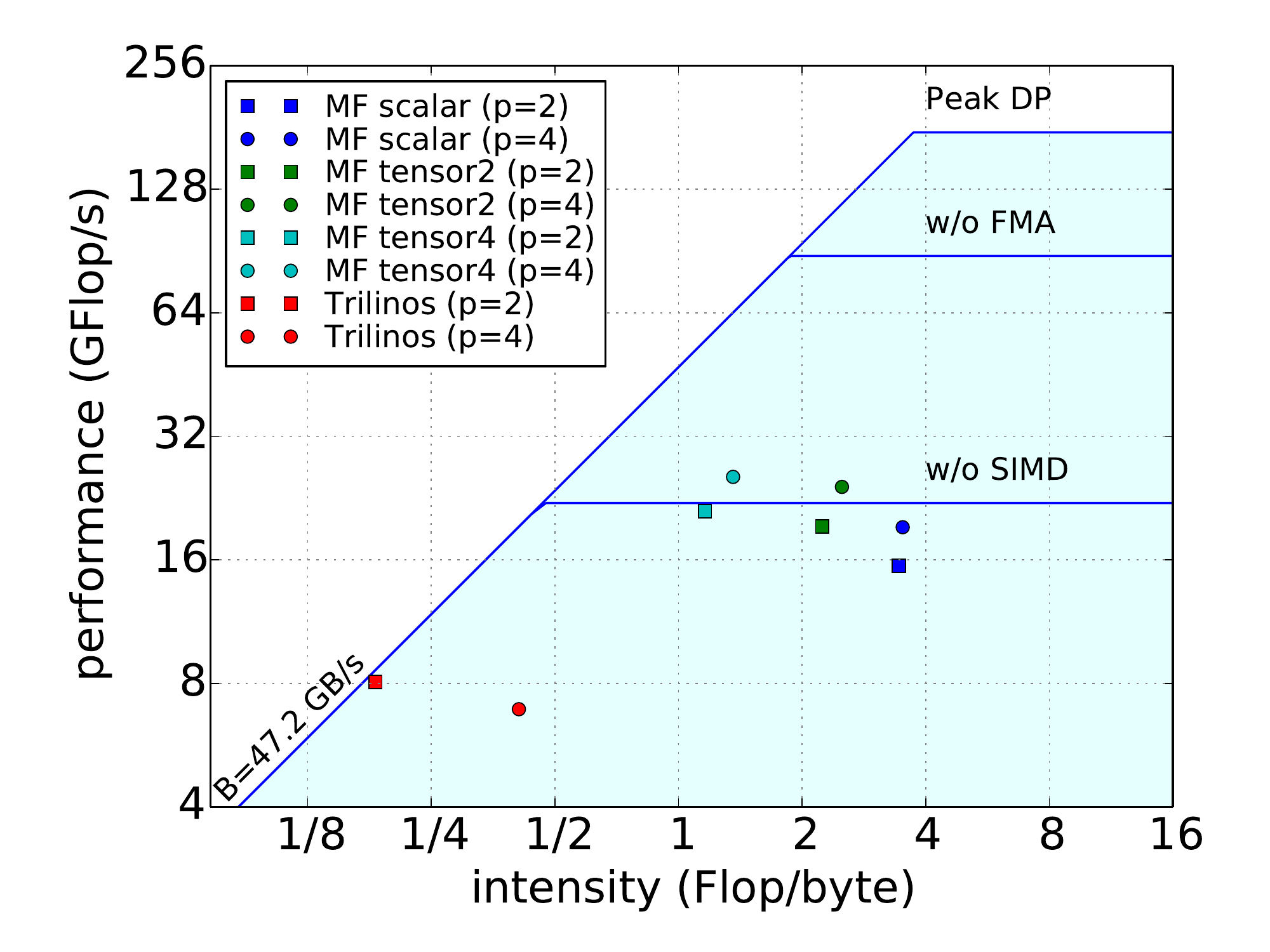}
    \caption{3D}
    \label{fig:roofline_3d}
  \end{subfigure}
  \caption{Roofline model for a selection of polynomial degrees. `Peak DP' denotes the peak double precision performance ($P$),
  `w/o FMA' represents a ceiling without Fused Multiply-Add ($P/2$), finally `w/o SIMD' is a ceiling where both FMA
    and SIMD vectorization are discarded ($P/8$).}%
  \label{fig:roofline}
\end{figure}

Figure \ref{fig:roofline} shows results of the roofline performance model.
As expected, the matrix-based approach is memory bound and (on average) achieves the performance of $\pgfmathprintnumber{8.2093245}\,\rm{GFlop/s}$ (4.6\% of the peak performance) in 2D.
The matrix-free implementations of the fine-strain elastic tangent operator (on average) lead to a much better performance of
$\pgfmathprintnumber{18.8047083222}\,\rm{GFlop/s}$ in 2D, as well as about an order of magnitude higher computational intensity, but still stays in the memory-bound regime.
The achieved performance is about 10.6\% of the peak arithmetic performance. This is much lower than the
70\% reported for the Laplace operator on Cartesian meshes in \cite{kronbichler12}, but similar to what was obtained in \cite[Figure 14]{kronbichler2017fast} when not optimizing for Cartesian meshes.
Unfortunately, performance measurements for non-Cartesian meshes were not reported in \cite{kronbichler12}.

When comparing the data for ``tensor4'' caching strategy (Algorithm \ref{alg:mf_tensor4}) to that from the ``tensor2'' (Algorithm \ref{alg:mf_tensor2}) caching strategy,
we observe a marginal increase in computational intensity and almost unchanged performance.
This coincides with our expectation that the latter algorithm requires transferring fewer bytes to perform operations on a single quadrature point
(compare line 9 in Algorithm \ref{alg:mf_tensor4} that loads a symmetric rank-4 tensor to the line 9 in Algorithm \ref{alg:mf_tensor2} that evaluates the
same contraction by loading two scalars $c_1$ and $c_2$).
On the other hand, the ``scalar'' caching strategy (Algorithm \ref{alg:mf_scalar}) demonstrates a higher computational intensity.
This is not surprising as at each quadrature point we need to re-evaluate the Kirchoff stress $\gz \tau$.
Given that in this case we need to evaluate gradients with respect to both deformed and undeformed configurations, an additional
second-order tensor per quadrature point (inverse of the Jacobian transformation) will be cached as compared to the Algorithm \ref{alg:mf_tensor2}.
As a result both algorithms will have comparable memory requirements (also confirmed by Figures \ref{fig:benchmark_miehe_Emmy_memory2}, \ref{fig:benchmark_miehe_Emmy_memory3}, \ref{fig:benchmark_miehe_IWR_memory2} and \ref{fig:benchmark_miehe_IWR_memory3}),
however the ``scalar'' algorithm performs a large number of repetitive operations that are
independent on the right hand side vector (steps 8--11 in Algorithm \ref{alg:mf_scalar}), which makes it the worst option among the three.

\begin{figure}[!ht]
  \centering
  \begin{subfigure}[b]{0.49\textwidth}
      \centering
      \includegraphics[width=\textwidth]{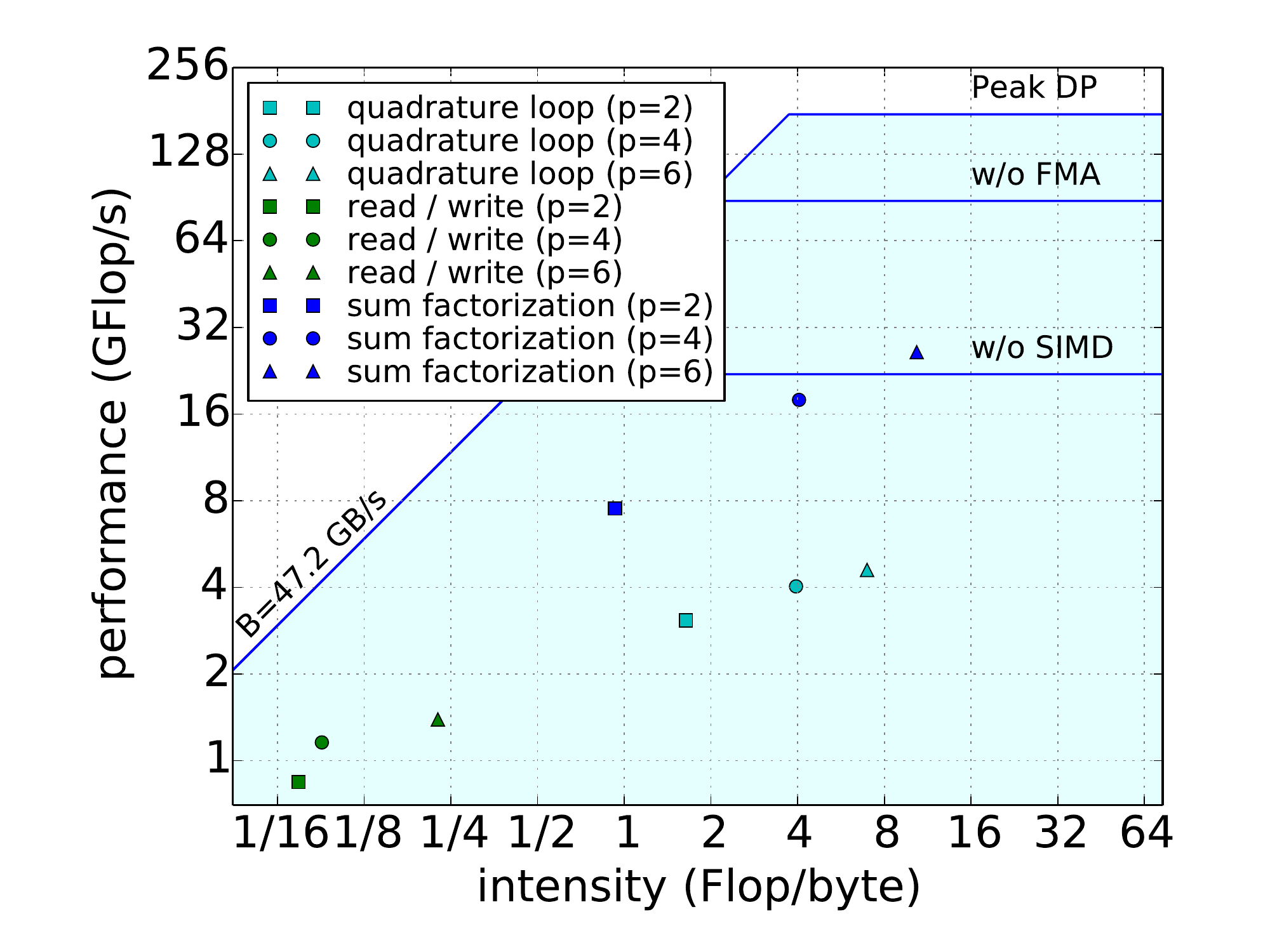}
      \caption{Roofline (2D)}
      \label{fig:roofline_2d_tensor4}
  \end{subfigure}
  \begin{subfigure}[b]{0.49\textwidth}
    \centering
    \includegraphics[width=\textwidth]{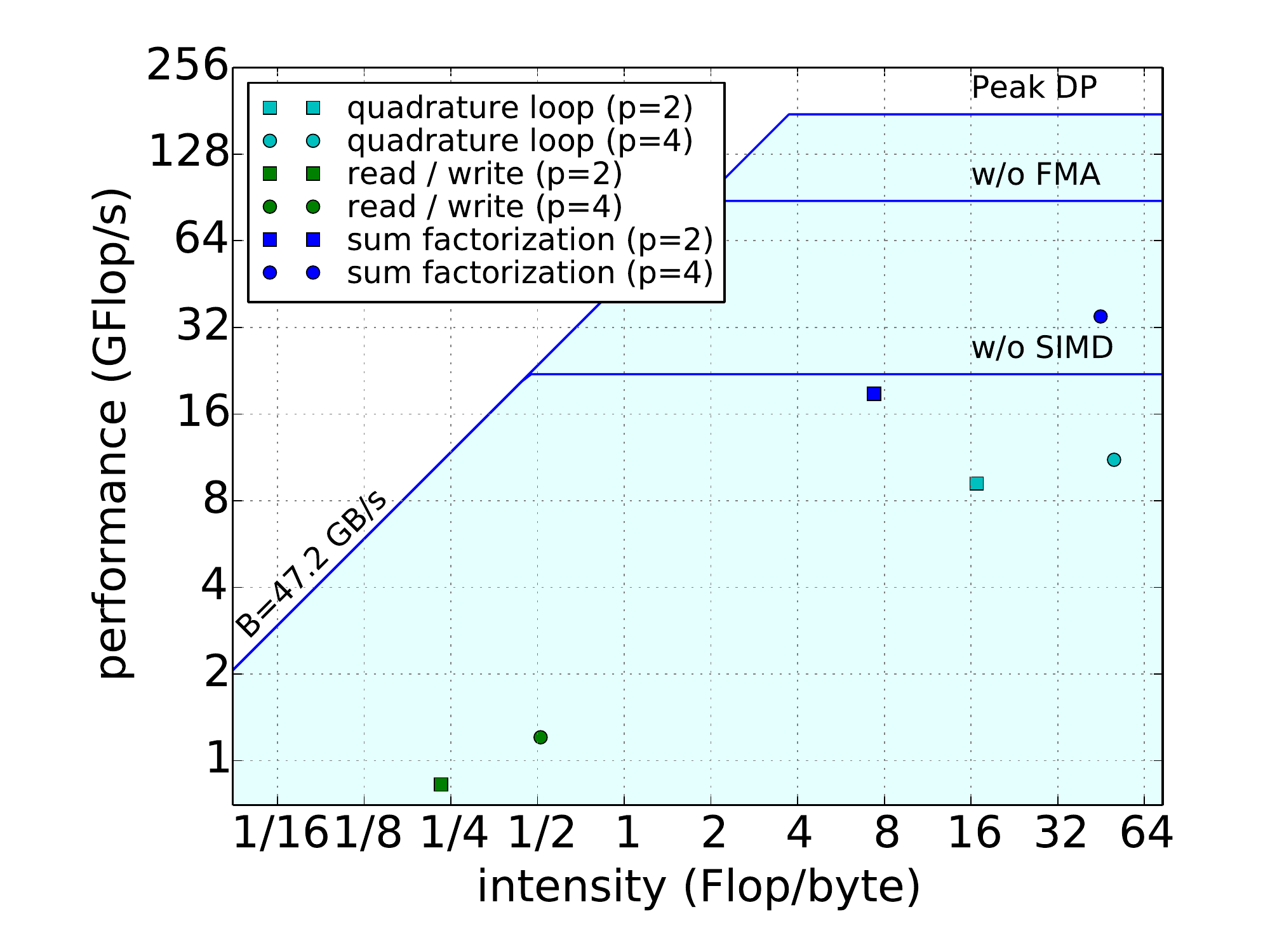}
    \caption{Roofline (3D)}
    \label{fig:roofline_3d_tensor4}
  \end{subfigure}
  ~
  \begin{subfigure}[b]{0.49\textwidth}
      \centering
      \includegraphics[width=\textwidth]{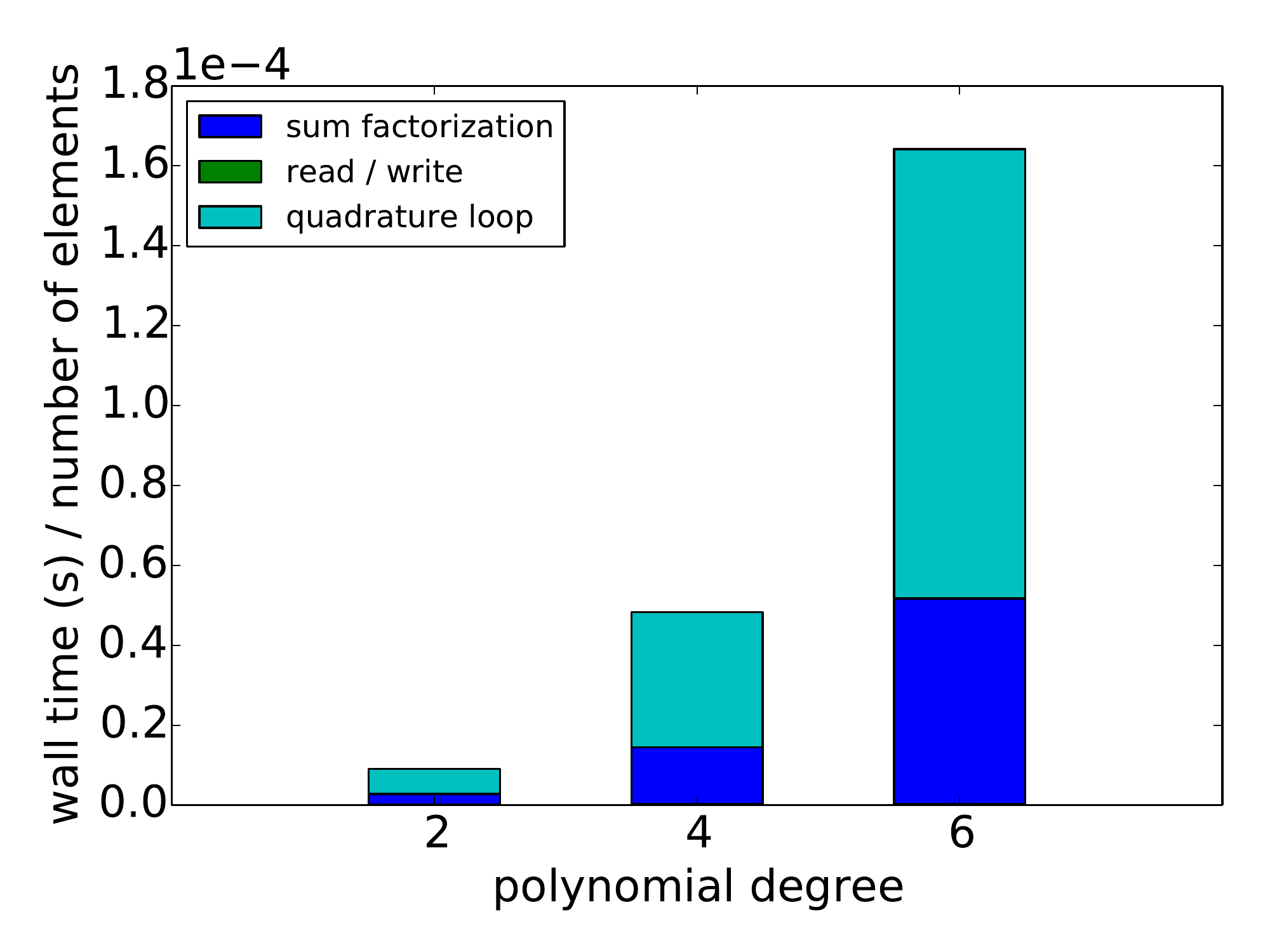}
      \caption{Computing times (2D)}
      \label{fig:breakdown_stackedbar_2d_tensor4}
  \end{subfigure}
  \begin{subfigure}[b]{0.49\textwidth}
    \centering
    \includegraphics[width=\textwidth]{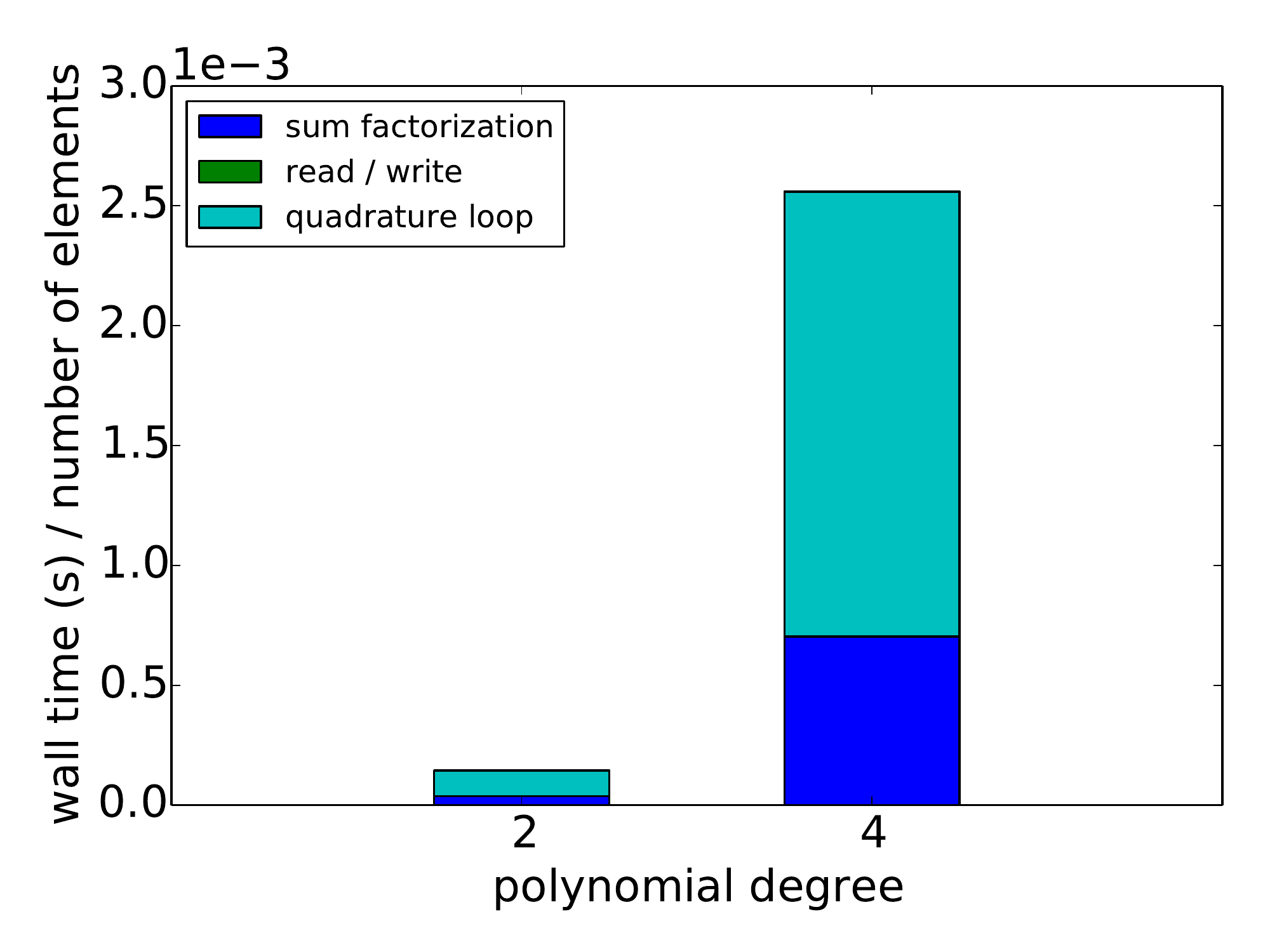}
    \caption{Computing times (3D)}
    \label{fig:breakdown_stackedbar_3d_tensor4}
  \end{subfigure}
  \caption{Breakdown of computing times and roofline analysis of various steps of Algorithm \ref{alg:mf_tensor4}.}%
  \label{fig:breakdown}
\end{figure}

In order to have a better understanding of the behavior of the most promising ``tensor4'' implementation,
we collect its performance data and computing times for various stages inside the element loop (see Algorithm \ref{alg:mf_tensor4}):
(i) `read / write' denotes reading and writing from/to a global vector (lines 2 and 13);
(ii) `sum factorization' denotes application of sum factorization techniques (lines 4, 5 and 12);
(iii) `quadrature loop' denotes operations performed on each quadrature (lines \mbox{6 -- 11}).
Given that the measurements are now down within the cell loop, we collected results on smaller meshes by adopting only one
global refinement in 2D and no global mesh refinement in 3D.
It is clear from Figure \ref{fig:breakdown_stackedbar_2d_tensor4} and  \ref{fig:breakdown_stackedbar_3d_tensor4} that the
`read / write' steps take a negligible amount of time on each element.
For higher order bases, the majority of time is spent within the quadrature loop.
When analyzing the roofline data (Figure \ref{fig:roofline_2d_tensor4} and \ref{fig:roofline_3d_tensor4}),
we observe that this part of the algorithm demonstrates rather low performance, which in our opinion is the source of the
suboptimal performance of the whole matrix-free matrix-vector product (see Figure \ref{fig:roofline}).
In this case the, quadrature loop consists of, essentially, a
single contraction between two second order tensors (step 7 in Algorithm \ref{alg:mf_tensor4}) and a
double contraction between a symmetric second order tensor and a fourth order tensor (step 9 in Algorithm \ref{alg:mf_tensor4}).

\begin{table}
  \centering
  \resizebox{\textwidth}{!}{
  \begin{tabular}{|c|cc|ccc|ccc|ccc|}
  \hline
                & \multicolumn{2}{c|}{serial} & \multicolumn{3}{c|}{MPI} & \multicolumn{3}{c|}{SIMD} & \multicolumn{3}{c|}{MPI+SIMD}  \\
  \hline
  p             & time  & GFlop/s              & time & GFlop/s & speedup & time & GFlop/s & speedup & time & GFlop/s & speedup \\
  \hline
  2& \pgfmathprintnumber{0.71685} & \pgfmathprintnumber{0.8747712} & \pgfmathprintnumber{0.07427} & \pgfmathprintnumber{8.57318} & \pgfmathprintnumber{9.65194560388} & \pgfmathprintnumber{0.3939} & \pgfmathprintnumber{1.6528706} & \pgfmathprintnumber{1.81987814166} & \pgfmathprintnumber{0.04142} & \pgfmathprintnumber{16.0073657} & \pgfmathprintnumber{17.306856591} \\
  4& \pgfmathprintnumber{0.47603} & \pgfmathprintnumber{1.3088965} & \pgfmathprintnumber{0.0511} & \pgfmathprintnumber{12.3860771} & \pgfmathprintnumber{9.3156555773} & \pgfmathprintnumber{0.25747} & \pgfmathprintnumber{2.3690572} & \pgfmathprintnumber{1.84887559716} & \pgfmathprintnumber{0.02786} & \pgfmathprintnumber{22.2420699} & \pgfmathprintnumber{17.0865039483} \\
  6& \pgfmathprintnumber{0.24683} & \pgfmathprintnumber{1.5497601} & \pgfmathprintnumber{0.02636} & \pgfmathprintnumber{14.8125524} & \pgfmathprintnumber{9.36380880121} & \pgfmathprintnumber{0.13089} & \pgfmathprintnumber{2.8544567} & \pgfmathprintnumber{1.88578195431} & \pgfmathprintnumber{0.01487} & \pgfmathprintnumber{26.0509717} & \pgfmathprintnumber{16.5991930061} \\
  \hline
  \end{tabular}
  }
  \caption{Wall-clock time in seconds and performance in GFlops of Algorithm \ref{alg:mf_tensor4} in 2D for various combinations of polynomial degrees,
  vectorization and parallelization.}
  \label{tab:numbers_2d}
\end{table}

\begin{table}
  \centering
  \resizebox{\textwidth}{!}{
  \begin{tabular}{|c|cc|ccc|ccc|ccc|}
  \hline
                & \multicolumn{2}{c|}{serial} & \multicolumn{3}{c|}{MPI} & \multicolumn{3}{c|}{SIMD} & \multicolumn{3}{c|}{MPI+SIMD}  \\
  \hline
  p             & time  & GFlop/s              & time & GFlop/s & speedup & time & GFlop/s & speedup & time & GFlop/s & speedup \\
  \hline
  2& \pgfmathprintnumber{1.82927} & \pgfmathprintnumber{1.4408182} & \pgfmathprintnumber{0.19184} & \pgfmathprintnumber{13.75528} & \pgfmathprintnumber{9.5353940784} & \pgfmathprintnumber{1.24738} & \pgfmathprintnumber{2.2390063} & \pgfmathprintnumber{1.46648976254} & \pgfmathprintnumber{0.13397} & \pgfmathprintnumber{21.0187594} & \pgfmathprintnumber{13.6543255953} \\
  4& \pgfmathprintnumber{1.14988} & \pgfmathprintnumber{1.7483867} & \pgfmathprintnumber{0.1211} & \pgfmathprintnumber{16.6162812} & \pgfmathprintnumber{9.49529314616} & \pgfmathprintnumber{0.74974} & \pgfmathprintnumber{2.7527409} & \pgfmathprintnumber{1.53370501774} & \pgfmathprintnumber{0.08166} & \pgfmathprintnumber{25.4865287} & \pgfmathprintnumber{14.0813127602} \\
  \hline
  \end{tabular}
  }
  \caption{Wall-clock time in seconds and performance in GFlops of Algorithm \ref{alg:mf_tensor4} in 3D for various combinations of polynomial degrees,
  vectorization and parallelization.}
  \label{tab:numbers_3d}
\end{table}

Table \ref{tab:numbers_2d} and Table \ref{tab:numbers_3d} show the speed-up obtained due to explicit SIMD vectorization intrinsics and due to MPI. For these rather small count of MPI processes, the results confirm that the performance gained is independent. Furthermore, the improvement due to explicit intrinsics for the data used in the quadrature loop, is comparable to the ones in \cite[Figure 6]{kronbichler2017fastcomputer} where the speedup factor is $1.2-1.8$ for a similar architecture.

To summarize this section, we can attribute the speed-up of the matrix-free finite-strain tangent operators for higher-order elements to two factors.
The first one is the higher algorithmic intensity as compared to the matrix-based implementation, which is clearly memory
bound, see Figure \ref{fig:roofline}.
The second one is the reduction in the number of floating point operations per DoF thanks to the sum factorization technique, as discussed in
Section \ref{sec:mf} and \cite[Section 2.4]{kronbichler12}.
Among the three suggested implementations, the Algorithm \ref{alg:mf_tensor4} is the most flexible approach with relatively high computational intensity.

\subsection{Preconditioned iterative solver}

\begin{figure}[!ht]
  \begin{subfigure}[b]{0.49\textwidth}
    \centering
    \includegraphics[width=\textwidth]{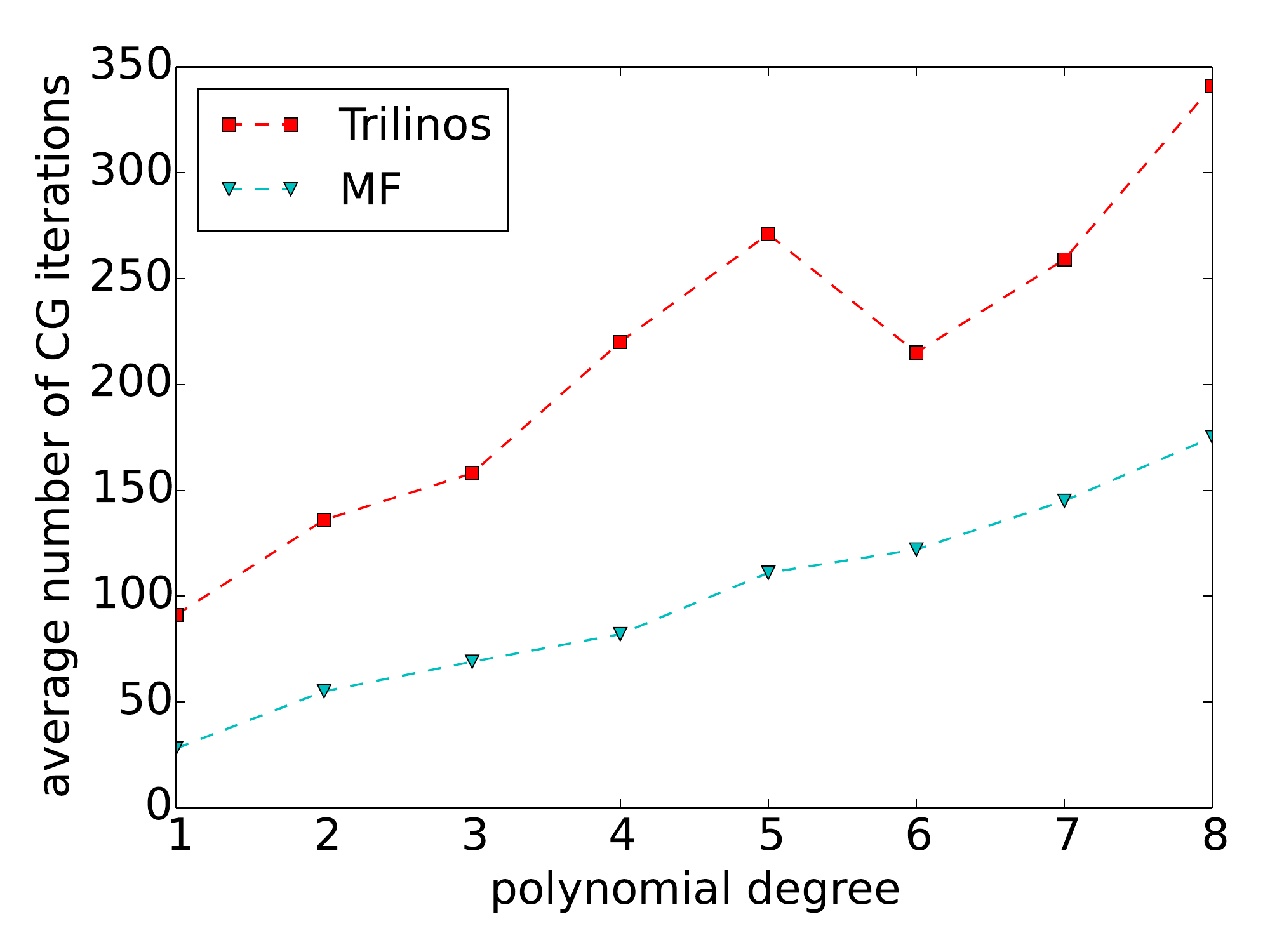}
    \caption{CG iterations (2D)}
    \label{fig:benchmark_miehe_Emmy_cg2}
  \end{subfigure}
  \begin{subfigure}[b]{0.49\textwidth}
    \centering
    \includegraphics[width=\textwidth]{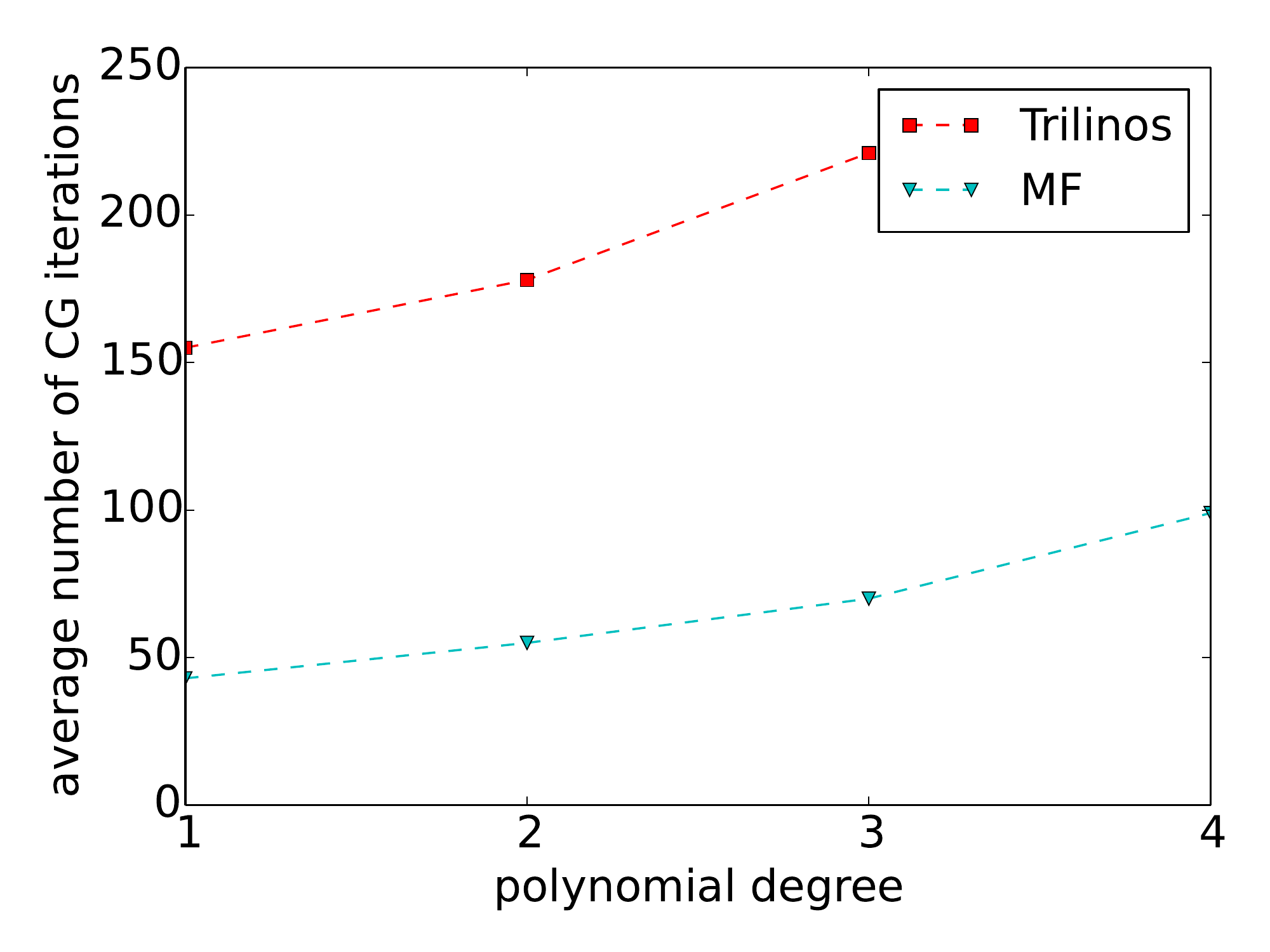}
    \caption{CG iterations (3D)}
    \label{fig:benchmark_miehe_Emmy_cg3}
  \end{subfigure}
  ~
  \begin{subfigure}[b]{0.49\textwidth}
    \centering
    \includegraphics[width=\textwidth]{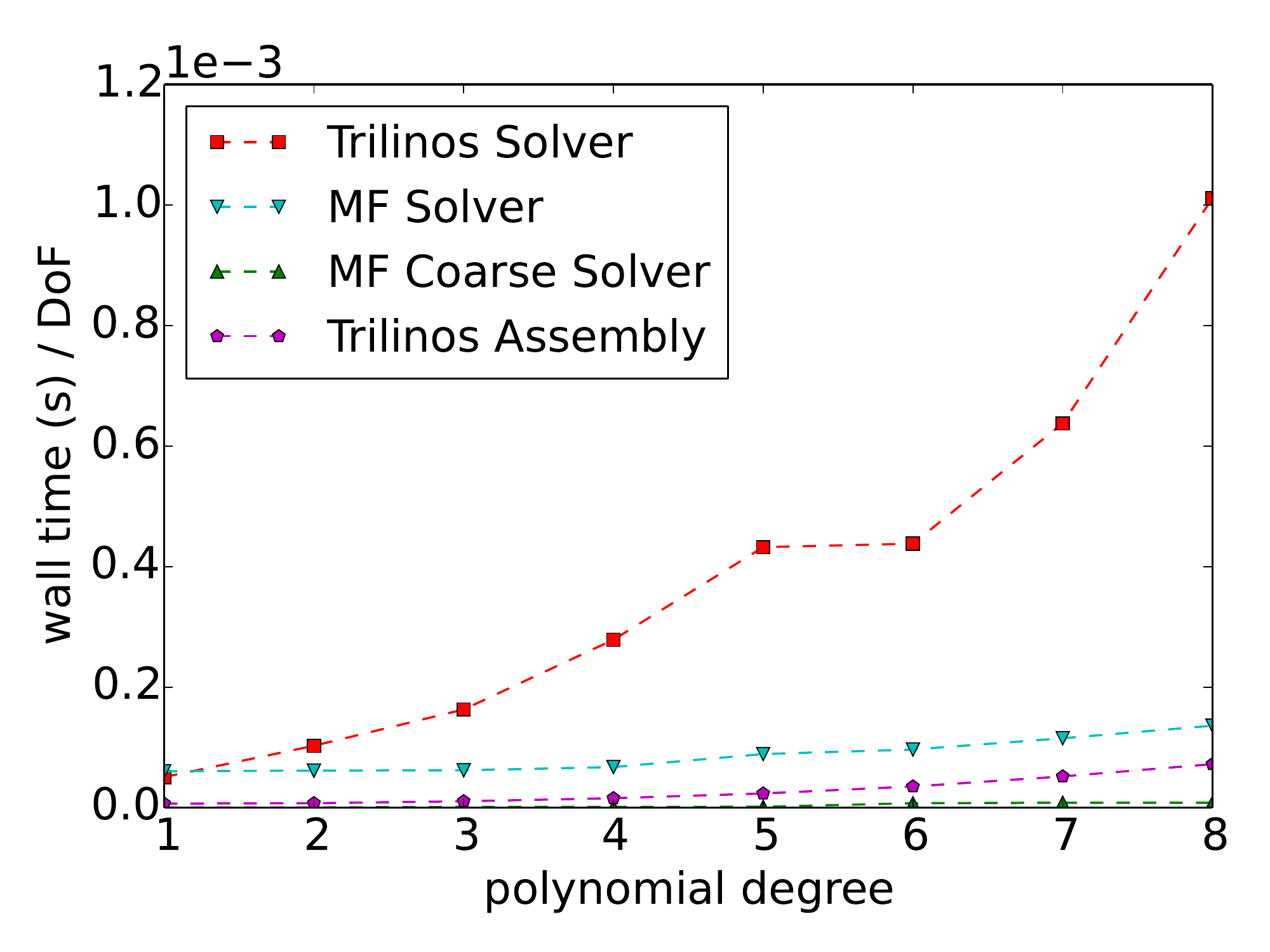}
    \caption{CG solution time (2D)}
    \label{fig:benchmark_miehe_Emmy_sol2}
  \end{subfigure}
  \begin{subfigure}[b]{0.49\textwidth}
    \centering
    \includegraphics[width=\textwidth]{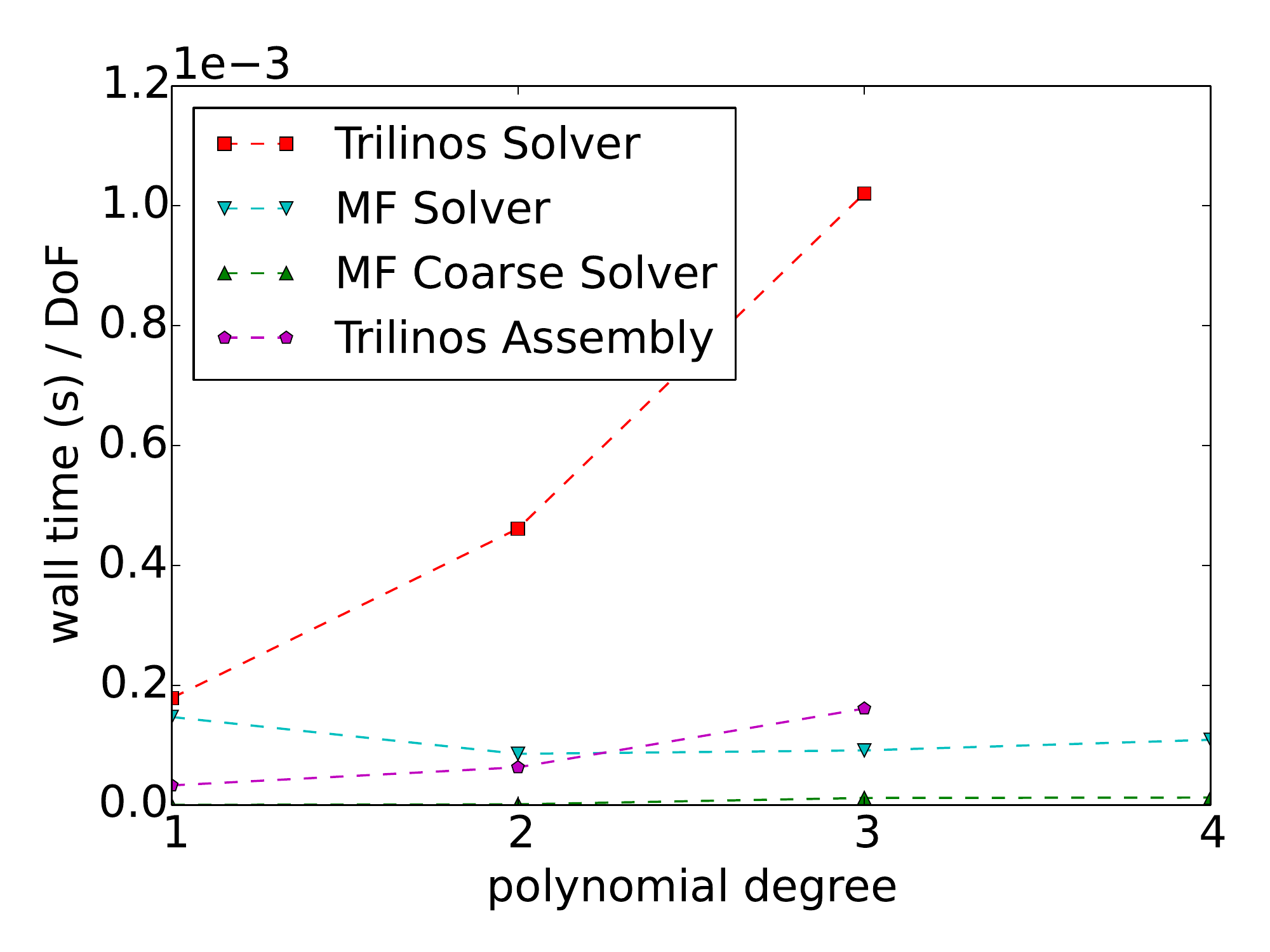}
    \caption{CG solution time (3D)}
    \label{fig:benchmark_miehe_Emmy_sol3}
  \end{subfigure}
  \caption{Emmy cluster, RRZE, iterative solver.}%
  \label{fig:benchmark_miehe_Emmy_cg}
\end{figure}

\begin{figure}[!ht]
  \begin{subfigure}[b]{0.49\textwidth}
    \centering
    \includegraphics[width=\textwidth]{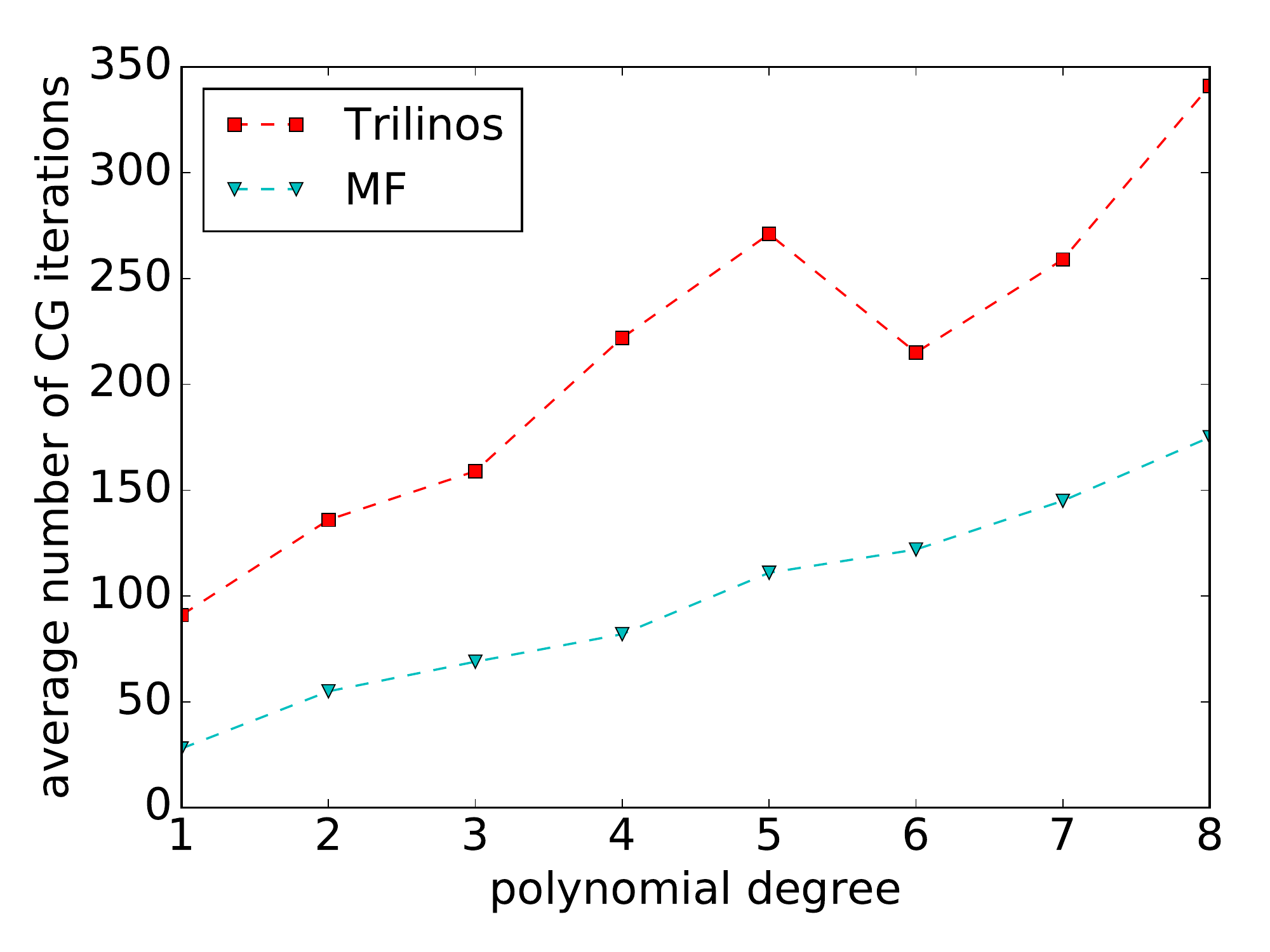}
    \caption{CG iterations (2D)}
    \label{fig:benchmark_miehe_IWR_cg2}
  \end{subfigure}
  \begin{subfigure}[b]{0.49\textwidth}
    \centering
    \includegraphics[width=\textwidth]{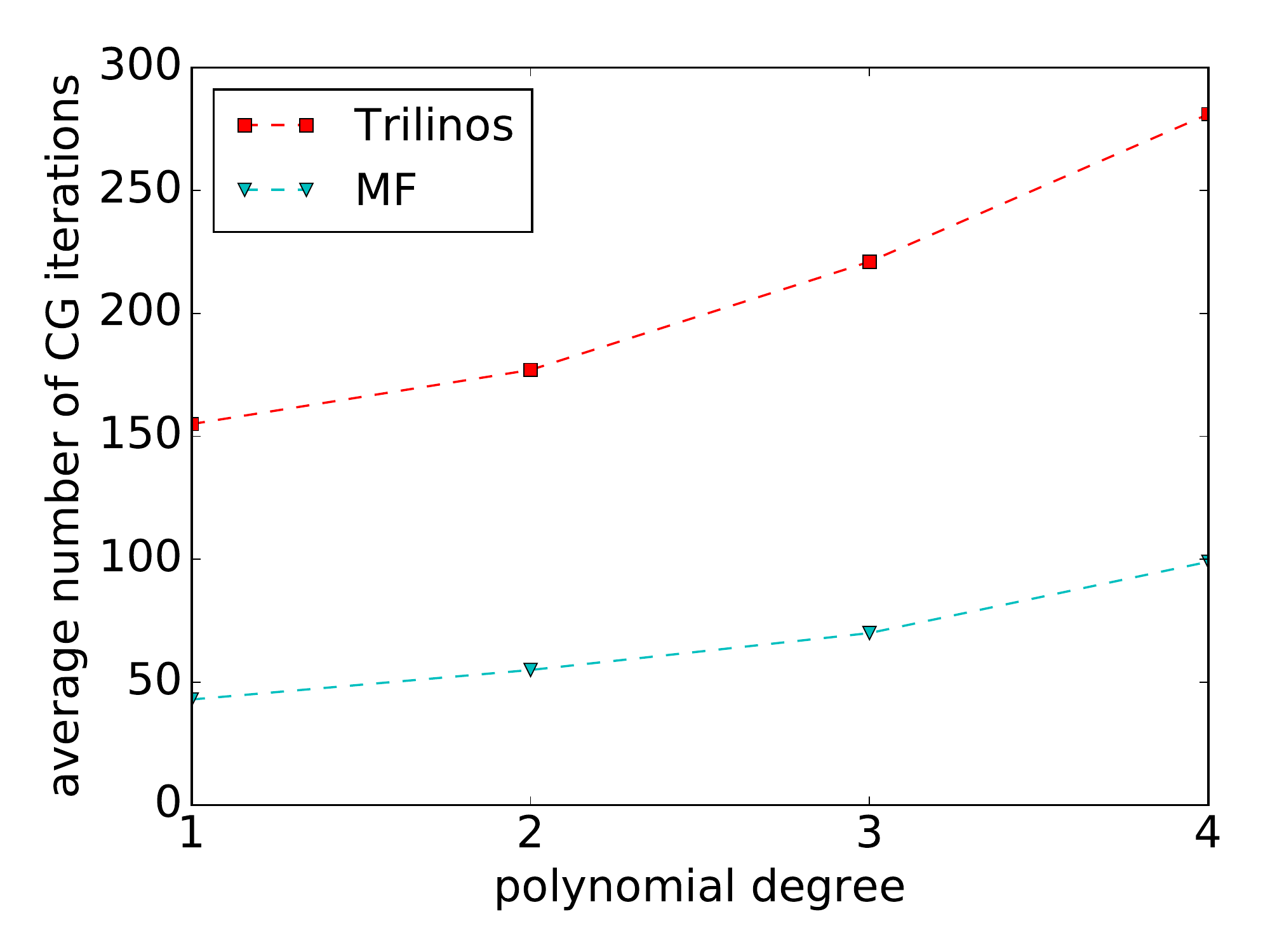}
    \caption{CG iterations (3D)}
    \label{fig:benchmark_miehe_IWR_cg3}
  \end{subfigure}
  ~
  \begin{subfigure}[b]{0.49\textwidth}
    \centering
    \includegraphics[width=\textwidth]{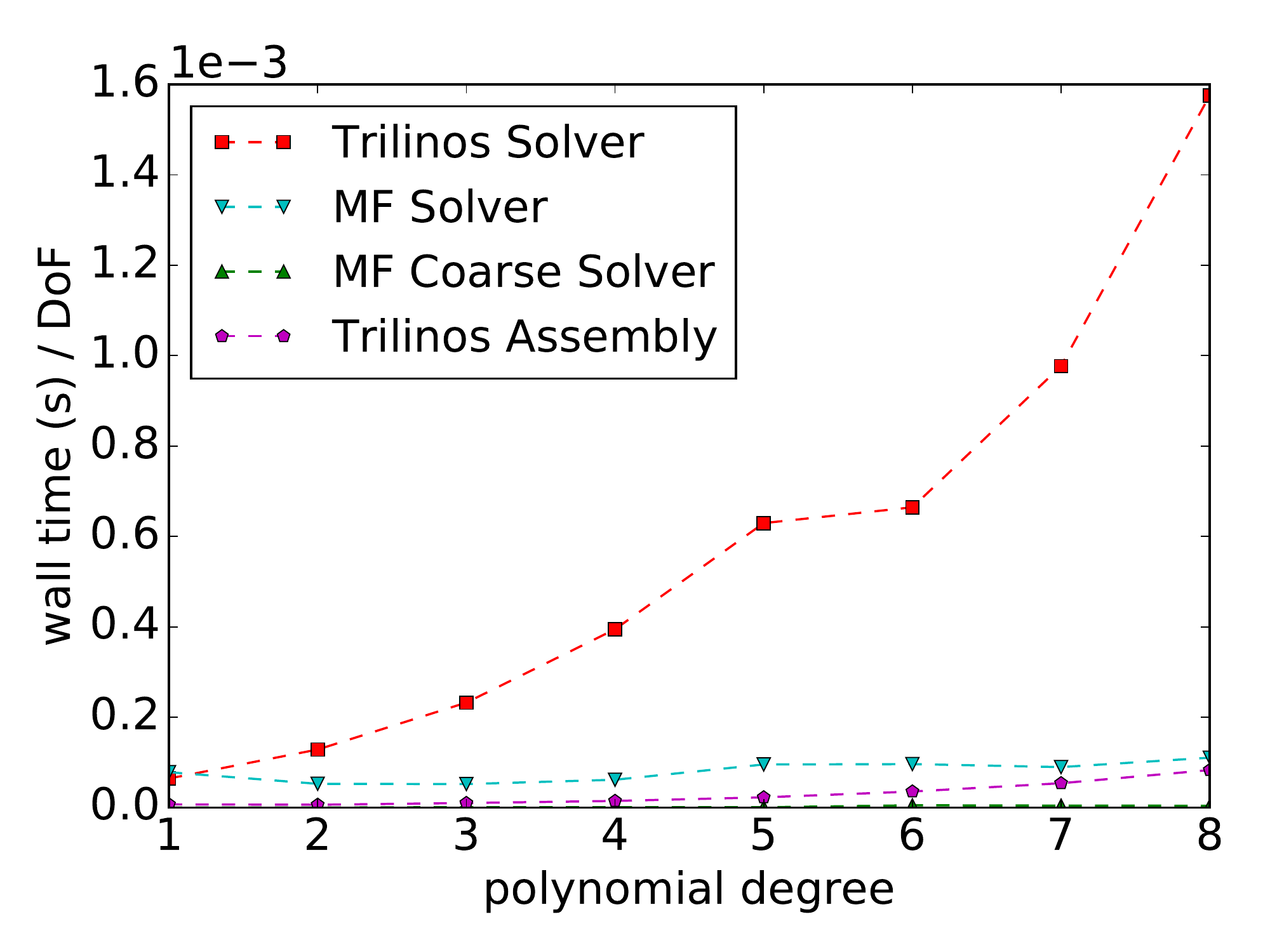}
    \caption{CG solution time (2D)}
    \label{fig:benchmark_miehe_IWR_sol2}
  \end{subfigure}
  \begin{subfigure}[b]{0.49\textwidth}
    \centering
    \includegraphics[width=\textwidth]{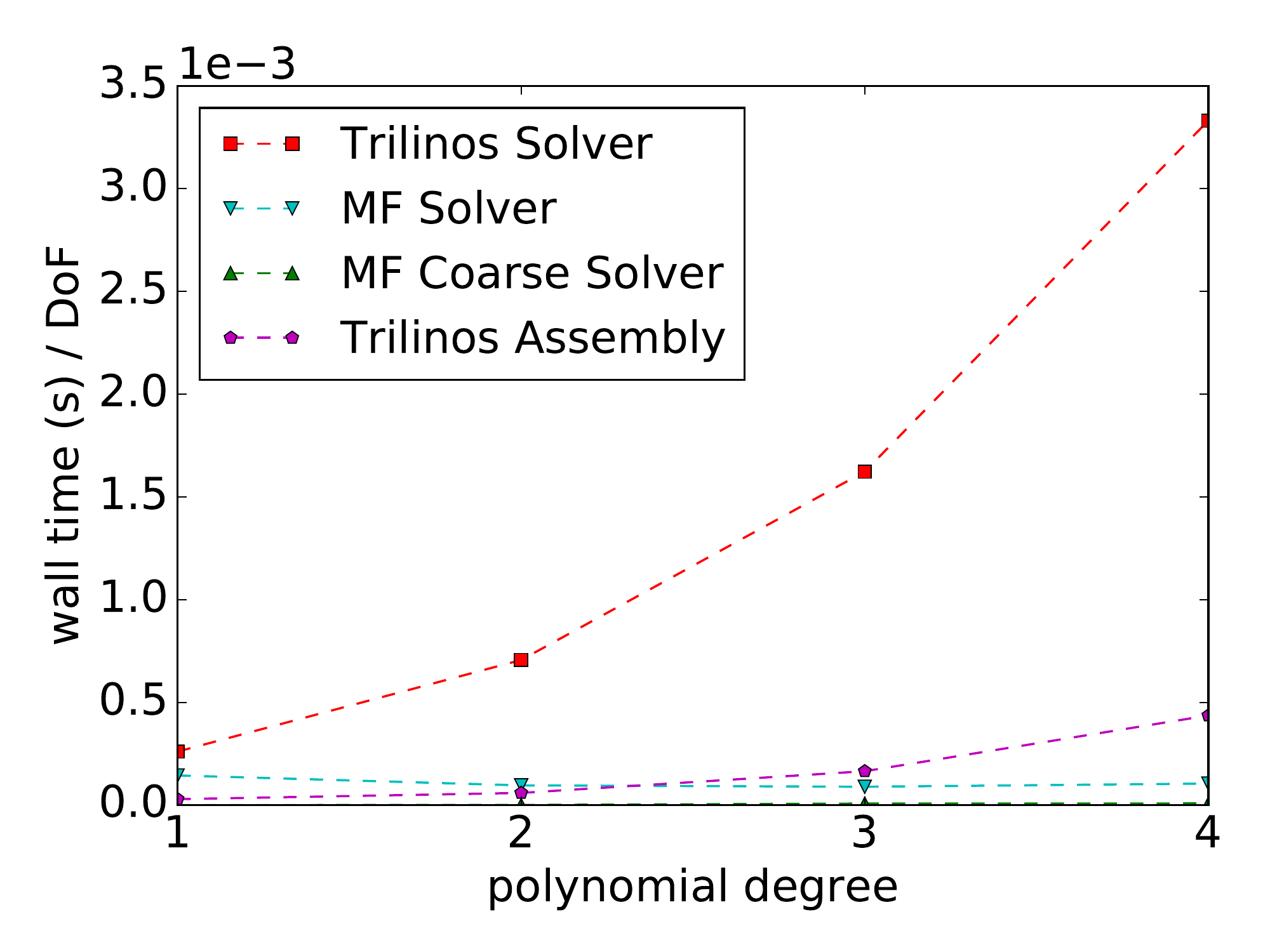}
    \caption{CG solution time (3D)}
    \label{fig:benchmark_miehe_IWR_sol3}
  \end{subfigure}
  \caption{IWR cluster.}%
  \label{fig:benchmark_miehe_IWR_cg}
\end{figure}

Next, we evaluate the performance of the proposed geometric multigrid preconditioner.
Our main goal is to examine whether or not the adopted level operators lead to an efficient preconditioner from the linear algebra perspective. To that end,
we consider the average number of CG iterations throughout the entire simulation (i.e. each Newton-Raphson iteration and each loading step).
Figures \ref{fig:benchmark_miehe_Emmy_cg2}, \ref{fig:benchmark_miehe_Emmy_cg3}, \ref{fig:benchmark_miehe_IWR_cg2} and \ref{fig:benchmark_miehe_IWR_cg3} show a
 noticeable increase in the average number of CG iterations for different polynomial degrees
 for the GMG preconditioner.
Although not reported in this section, we identified that
the coarse level solver setup will greatly affect the performance of the preconditioner.
In particular, we can improve the number of iterations by choosing a more accurate coarse level solver,
however this leads to an overall more expensive linear solver.
However, the settings used in this study are optimized for the considered problem in terms of the wall-time spent in the linear solver.
The black-box AMG preconditioner
demonstrates the same trend in terms of the number of CG iterations, but
requires many more iterations for convergence.
Consequently, we can conclude that for the here considered finite strain elasticity problem (i) the selected smoother is suitable, and that (ii) although the adopted level operators are not built using the triple product $\gz A^{l+1}=\gz I^{l+1}_{l} \gz A^l \gz I^l_{l+1}$, they still result in a multigrid preconditioner which requires much fewer CG iterations as compared to the algebraic multigrid.

Additionally, we compare solution times of the preconditioned iterative solvers.
Our goal is to determine whether the matrix-free approach with a GMG preconditioner can be faster than the matrix-based approach with a well-established implementation of an AMG preconditioner.
We note specifically that the comparison is restricted to the choice of level operators used with the pre-existing GMG framework, and not the framework itself.
Figures \ref{fig:benchmark_miehe_Emmy_sol2}, \ref{fig:benchmark_miehe_Emmy_sol3}, \ref{fig:benchmark_miehe_IWR_sol2} and \ref{fig:benchmark_miehe_IWR_sol3} confirm that the efficiency of the proposed GMG preconditioner translates into faster solution times for all polynomial degrees.
Most importantly, the wall time per DoF for the matrix-free solver with GMG preconditioner stays almost constant both for 2D and 3D problems.
In comparison, the wall time of the matrix-based iterative solver grows rapidly with increasing polynomial degree.
Clearly this is related both to the performance of the preconditioner from the linear algebra perspective as well as the matrix-vector products, studied in the previous section.
Based on these results, we conclude that the matrix-free approach with the adopted here geometric multigrid preconditioner can deliver a very competitive solution strategy for engineering problems
with compressible hyperelastic finite strain material models.

Extension of this study from the node level to the cluster level is beyond the scope of this work.
However, as a step in this direction we have performed a preliminary study of the weak scalability of the implementation.
Table \ref{tab:weak_3d} demonstrates a good weak scaling up to 1280 processes\footnote{%
This corresponds to 64 nodes, which is the maximum we could request on the RRZE Emmy HPC cluster.
} of the wall time per matrix-vector product and per iteration in the CG solver for the 3D problem with a tri-quadratic basis.
This confirms that the matrix-free and multigrid frameworks implemented in the \texttt{deal.II} library can provide a competitive solution for cluster-sized problems as well,
also compare \cite[Figure 7]{Krank2017} for scaling results for up to 147,456 CPU cores and 34.4 billion degrees of freedom for an incompressible flow problem.

\begin{table}
  \centering
  \begin{tabular}{|r|r|c|r|c|}
  \hline
  cores  & $N_{DoF}$ & $N_{gref}$ & vmult [s] & CG [s / iteration] \\
  \hline
  20 & 4,442,880 & 3 & \pgfmathprintnumber{0.0675} & \pgfmathprintnumber{0.529553679131} \\
  160 & 35,071,488 & 4 & \pgfmathprintnumber{0.0705} & \pgfmathprintnumber{0.562637362637} \\
  1280 & 278,694,912 & 5 & \pgfmathprintnumber{0.083} & \pgfmathprintnumber{0.688792165397} \\
  \hline
  \end{tabular}
  \caption{Weak scaling of Algorithm \ref{alg:mf_tensor4} in 3D for quadratic polynomial basis.}
  \label{tab:weak_3d}
\end{table}

\section{Summary and Conclusions}
\label{sec:summary}

In this contribution, we proposed and numerically investigated three different matrix-free implementations of tangent operators for finite-strain elasticity with heterogeneous materials.
To the best of our knowledge, this is the first work that applies matrix-free sum factorization techniques to partial differential equations arising in this case.
Among the three examined algorithms,
the implementation that caches the material part of the fourth-order spatial tangent stiffness tensor together with the second-order Kirchhoff stress was shown to be faster than the matrix-based approach for polynomial degrees higher than one in 3D and polynomial degrees higher than two in 2D.
Given that we did not utilize the specific form of the material part of the stiffness tensor in this case, we conclude
that the matrix-free implementation of the tangent operator can be applied to any constitutive model which operates with the material part of the fourth-order spatial tangent stiffness tensor on the quadrature point level.

The roofline model indicates that the matrix-free finite-strain elasticity tangent operators have an order magnitude higher algorithmic intensity.
We identified that one of the bottlenecks is related to the quadrature loop on each element where double contractions of fourth-order and second-order symmetric tensor
as well as single contraction of two second-order tensors is performed.

We also propose a method by which to construct level tangent operators and employ them to define a geometric multigrid preconditioner
with standard geometric transfer operations between each level.
The GMG was applied to heterogeneous material assuming that the coarsest level can provide an adequate discretization of the heterogeneity.
The numerical studies indicate that the proposed preconditioner
leads to much fewer iterations of the iterative solver as compared to the algebraic multigrid preconditioner.
The multigrid matrix-free preconditioner also
leads to a solution approach that is faster than the matrix-based AMG for all polynomial degrees studied here.
Most importantly, the wall time per degree of freedom for solving the problem is close to being constant.
On the other hand, the matrix-based implementation with AMG becomes prohibitively expensive for higher-order bases.
We conclude that the matrix-free implementations of tangent operators for finite-strain elasticity together with
the geometric multigrid preconditioner is a very competitive solution strategy, as compared to more traditional matrix-based approach.
Our future work will be focused on adopting the proposed matrix-free multigrid solution approach to $\operatorname{FE}^2$ homogenization
as well studying its behavior on the cluster level.

\ifijnme
\acks
\else
\section*{Acknowledgements}
\fi

D.~Davydov acknowledges the financial support of the German Research Foundation (DFG), grant DA 1664/2-1.

D.~Davydov, D.~Arndt and J-P.~Pelteret are grateful to Martin Kronbichler (TU Munich), Jed Brown (CU Boulder) and Veselin Dobrev (Lawrence Livermore National Laboratory) for fruitful discussions on matrix-free operator evaluation approaches.

D.~Arndt was supported by the German Research Foundation (DFG) under the project ``High-order discontinuous
Galerkin for the exa-scale'' (\mbox{ExaDG}) within the priority program ``Software
for Exascale Computing'' (SPPEXA).

P.~Steinmann acknowledges the support of the Cluster of Excellence Engineering of Advanced Materials (EAM) which made this collaboration possible, as well as funding by the EPSRC Strategic Support Package ``Engineering of Active Materials by Multiscale/Multiphysics Computational Mechanics''.

\ifijnme
\bibliographystyle{wileyj}
\else
\fi


\end{document}